\documentclass[reqno,11pt]{amsart} 
\usepackage{amsmath}
\usepackage{amssymb}
\usepackage{amsthm}

\newcommand\NoBlackBoxes{\global\overfullrule0pt}
\NoBlackBoxes
\theoremstyle{plain} 

\def\4{\kern1pt}

\def\6{\vphantom0}

\def\8{\kern-10pt}
\def\7#1{_{(#1)}}

\begin{document}

\def\ffrac#1#2{\raise.5pt\hbox{\small$\4\displaystyle\frac{\,#1\,}{\,#2\,}\4$}}
\def\ovln#1{\,{\overline{\!#1}}}
\def\ve{\varepsilon}
\def\kar{\beta_r}

\title{R\'enyi divergence and the central limit theorem}

\author{S. G. Bobkov$^{1,4}$}
\thanks{1) School of Mathematics, University of Minnesota, USA;
Email: bobkov@math.umn.edu}
\address
{Sergey G. Bobkov \newline
School of Mathematics, University of Minnesota  \newline 
127 Vincent Hall, 206 Church St. S.E., Minneapolis, MN 55455 USA
\smallskip}
\email {bobkov@math.umn.edu} 

\author{G. P. Chistyakov$^{2,4}$}
\thanks{2) Faculty of Mathematics, University of Bielefeld, Germany;
Email: chistyak@math.uni-bielefeld.de}
\address
{Gennadiy P. Chistyakov\newline
Fakult\"at f\"ur Mathematik, Universit\"at Bielefeld\newline
Postfach 100131, 33501 Bielefeld, Germany}
\email {chistyak@math.uni-bielefeld.de}

\author{F. G\"otze$^{3,4}$}
\thanks{3) Faculty of Mathematics, University of Bielefeld, Germany;
Email: goetze@math.uni-bielefeld.de}
\thanks{4) Research partially supported by Humboldt Foundation, NSF grant
DMS-1612961 and SFB 701}
\address
{Friedrich G\"otze\newline
Fakult\"at f\"ur Mathematik, Universit\"at Bielefeld\newline
Postfach 100131, 33501 Bielefeld, Germany}
\email {goetze@mathematik.uni-bielefeld.de}


\begin{abstract}
We explore properties of the $\chi^2$ and more general R\'enyi (Tsallis) 
distances to the normal law. 
In particular we provide necessary and sufficient 
conditions for the convergence
to the normal law in the central limit theorem using these distances. 
Moreover,
 we derive exact rates of convergence in these distances  with respect to an increasing number of summands.
\end{abstract}

\maketitle
\markboth{S. G. Bobkov, G. P. Chistyakov and F. G\"otze}{
R\'enyi divergence from the normal law and the CLT}

\def\theequation{\thesection.\arabic{equation}}
\def\E{{\mathbb E}}
\def\R{{\mathbb R}}
\def\C{{\mathbb C}}
\def\P{{\mathbb P}}
\def\H{{\rm H}}
\def\Im{{\rm Im}}
\def\Tr{{\rm Tr}}

\def\k{{\kappa}}
\def\M{{\cal M}}
\def\Var{{\rm Var}}
\def\Ent{{\rm Ent}}
\def\O{{\rm Osc}_\mu}

\def\ep{\varepsilon}
\def\phi{\varphi}
\def\F{{\cal F}}
\def\L{{\cal L}}

\def\be{\begin{equation}}
\def\en{\end{equation}}
\def\bee{\begin{eqnarray*}}
\def\ene{\end{eqnarray*}}

\section{{\bf Introduction}}
\setcounter{equation}{0}

\vskip2mm
\noindent
Given random elements $X$ and $Z$ in a measurable space $(\Omega,\mu)$
with densities $p$ and $q$ (with respect to $\mu$), the $\chi^2$-distance of Pearson
$$
\chi^2(X,Z) = \int \frac{(p-q)^2}{q}\,d\mu
$$
represents an important measure of deviation of the distribution $P$ of $X$ from 
the distribution $Q$ of $Z$, which has been frequently used especally in Statistics 
and Information Theory (cf. e.g. [Le], [L-V], [V]). It is a rather strong 
distance-like quantity, which may be related to and included in the hierarchy 
of R\'enyi divergences (relative $\alpha$-entropies)
$$
D_\alpha(X||Z) = \frac{1}{\alpha - 1}\,\log \int
\Big(\frac{p}{q}\Big)^\alpha q\, d\mu \qquad (\alpha > 0)
$$
or equivalently, the R\'enyi divergence powers or the relative Tsallis entropies
$T_\alpha(X||Z)  = \frac{1}{\alpha - 1}\, [e^{(\alpha - 1) D_\alpha} - 1]$
(which do not depend on the choice of the dominating measure $\mu$).
The most important indexes are $\alpha = 0$, $\alpha = \frac{1}{2}$
(Hellinger distance), $\alpha = 1$ (Kullback-Leibler distance) and $\alpha = 2$
(quadratic R\'enyi/Tsallis divergence), in which case $T_2 = \chi^2$ and
$D_2 = \log(1 + \chi^2)$.

The functionals $D_\alpha$ and $T_\alpha$ are non-decreasing in $\alpha$, so, 
for growing indexes the distances are strengthening. In the range $0 < \alpha < 1$, 
all $D_\alpha$ are comparable to each other and are metrically equivalent 
to the total variation 
$\|P-Q\|_{\rm TV}$. However, the informational divergence $D = D_1 = T_1$ 
(called also entropic distance or relative entropy),
$$
D(X||Z) = \int p \log\frac{p}{q}\,d\mu,
$$
is much stronger, and this applies even more so to $D_\alpha$ with $\alpha>1$.
The difference between the different $D_\alpha$'s appears in applications like
the central limit theorem (CLT for short), which is studied in this paper. Here
we consider the $\chi^2$-divergence in the simplest situation of independent, 
identically distributed (i.i.d.) summands.

For i.i.d. random variables $X,X_1,X_2,\dots$ such that $\E X = 0$, $\E X^2 = 1$,
introduce the normalized sums
$$
Z_n = \frac{X_1 + \dots + X_n}{\sqrt{n}} \qquad (n = 1,2,\dots)
$$
together with their distributions $F_n$, which hence approach the standard 
normal law $\Phi$ in the weak sense. For convergence in the CLT using strong 
distances, recall that convergence in total variation was addressed in the 
1950's by Prokhorov [Pr]. 
He showed that $\|F_n - \Phi\|_{\rm TV}$ tends to zero as 
$n \rightarrow \infty$, if and only if $F_n$ has a non-trivial absolutely 
continuous component for some $n = n_0$, i.e., $\|F_{n_0} - \Phi\|_{\rm TV} < 2$ 
(in particular, this is true, if $X$ has density).
A similar description is due to Barron [B] in the 1980's for the 
Kullback-Leibler distance: $D(Z_n||Z)$ tends to zero for $Z \sim N(0,1)$, 
if and only if $D(Z_n||Z) < \infty$ for some $n = n_0$. The latter condition 
is fulfilled for a large family of underlying distributions, in particular, 
when $X$ has density $p$ such that
$$
\int_{-\infty}^\infty p(x) \log p(x)\,dx < \infty.
$$
Different aspects of such strong CLT's, including the non-i.i.d. situation and the
problem of rates or Berry-Esseen bounds, were studied by many
authors, and we refer an interested reader to 
[Li], [S-M], [A-B-B-N], [B-J], [J], [B-C-G2-4], [B-C-K], [B-C]. 

As for convergence in $D_\alpha$ with $\alpha > 1$, not much is known so far.
This case seems to be quite different in nature, and here the distance restricts 
the range of applicability of the CLT quite substantionally. When focusing 
on the particular value 
$\alpha = 2$, we are concerned with the behavior of the quantity
$$
\chi^2(Z_n,Z) = \int_{-\infty}^\infty \frac{(p_n(x) - \varphi(x))^2}{\varphi(x)}\,dx,
$$
where $p_n$ denotes the density of $Z_n$ and $\varphi$ is the standard normal
density. The finiteness of this integral already requires the existence of all 
moments of $X$ (and actually the existence of a ``Gaussian moment''). 
This condition is to be expected,
but the convergence to zero, and even the verification of the boundedness of
$\chi^2(Z_n,Z)$ in $n$ is rather delicate. This problem has been studied in the early 
1980's by Fomin [F] in terms of the exponential series (using Cramer's terminology)
for the density of $X$, 
$$
p(x) = \varphi(x) \sum_{k=1}^\infty \frac{\sigma_k}{2^k k!}\,H_{2k}(x),
$$
where $H_r$ denotes the $r$-th Chebyshev-Hermite polynomial. As a main result, 
he proved that $\chi^2(Z_n,Z) = O(\frac{1}{n})$ as $n \rightarrow \infty$, 
assuming that $p$ is compactly supported, symmetric, piecewise differentiable, 
such that the series coefficients satisfy $\sup_{k \geq 2} \sigma_k < 1$. 
This sufficient condition was verified for the uniform distribution on 
the interval $(-\sqrt{3},\sqrt{3})$ (this specific length is caused by 
the assumption $\E X^2 = 1$). However, for many other 
examples, Fomin's result does not seem to provide an applicable and
satisfactory answers.

Fortunately, more or less simple necessary and sufficient conditions can be 
stated for the convergence in $\chi^2$ by using the Laplace transform of 
the distribution of $X$. One of the purposes of this paper is to provide 
the following characterization of a class which may be called 
the ``domain of $\chi^2$-attraction to the normal law".

\vskip5mm
{\bf Theorem 1.1.} {\it We have $\chi^2(Z_n,Z)\to 0$ as $n\to\infty$, 
if and only if $\chi^2(Z_n,Z)$ is finite for some $n=n_0$, and
\be
\E\,e^{tX} < e^{t^2} \quad \text{for all real} \ t\ne 0.
\en
In this case the $\chi^2$-divergence admits an Edgeworth-type expansion
\be
\chi^2(Z_n,Z) \, = \,
\sum_{j=1}^{s-2} \frac{c_j}{n^j} + O\Big(\frac{1}{n^{s-1}}\Big) \quad as
\ n \to\infty,
\en
which is valid for every $s=3,4,\dots$ with coefficients $c_j$ representing 
certain polynomials in the moments $\alpha_k = \E X^k$, $k = 3,\dots,j+2$.
}

\vskip5mm
For $s=3$ this expansion simplifies to
$$
\chi^2(Z_n,Z) \, = \, \frac{\alpha_3^2}{6n} + O\Big(\frac{1}{n^2}\Big),
$$
and if $\alpha_3 = 0$ (as in the case of symmetric distributions), one may turn
to the next moment of order $s=4$, for which (1.2) yields
\be
\chi^2(Z_n,Z) \, = \, \frac{(\alpha_4 - 3)^2}{24\,n^2} + O\Big(\frac{1}{n^3}\Big).
\en

\vskip2mm
Let us note that the property $\chi^2(Z_n,Z) < \infty$ is rather close to 
the subgaussian condition (1.1). In particular, it implies that (1.1) is 
fulfilled for all $t$ large enough, as well as near zero due to the variance 
assumption. It may happen, however, that (1.1) is fulfilled for all $t \neq 0$ 
except just one value $t = t_0$ (and then there will be no CLT for the 
$\chi^2$-distance). Various examples illustrating these conditions together 
with the convergence in $\chi^2$ will be given in the end of the paper.

A similar characterization continues to hold in the {\it multidimensional} case
for mean zero i.i.d. random vectors $X,X_1,X_2,\dots$ in $\R^d$ normalized 
to have identity covariance. Here we endow the Euclidean space with the
canonical norm and scalar product. Moreover, one may extend these results to the
range of indexes $\alpha>1$, arriving at the following statement, where by
$\alpha^* = \frac{\alpha}{\alpha - 1}$ we denote the conjugate index.

\vskip5mm
{\bf Theorem 1.2.} {\it Let $Z$ denote a random vector in $\R^d$ having a standard 
normal distribution. Then $D_\alpha(Z_n||Z)\to 0$ as $n\to\infty$, 
if and only if $D_\alpha(Z_n||Z)$ is finite for some $n=n_0$, and
\be
\E\,e^{\left<t,X\right>} < e^{\alpha^* |t|^2/2} \quad \text{for all} \ \, 
t \in \R^d, \ t\ne 0.
\en
In this case, we necessarily have $D_\alpha(Z_n||Z) = O(1/n)$, and even 
$D_\alpha(Z_n||Z) = O(1/n^2)$, provided that the distribution of $X$ is 
symmetric about the origin.
}

\vskip5mm
Thanks to the existence of all moments of $X$, an Edgeworth-type expansion for 
$D_\alpha$ and $T_\alpha$ also holds similarly to (1.2), involving the mixed 
cumulants of the components of $X$. Such expansion shows in particular 
an equivalence
$$
D_\alpha(Z_n||Z) \sim T_\alpha(Z_n||Z) \sim \frac{\alpha}{2}\chi^2(Z_n,Z),
$$
provided that these distances tend to zero. Let us also note that 
the restriction imposed by (1.4) is asymptotically vanishing as $\alpha$ 
approaches 1. This means that we may expect to arrive at Barron's theorem 
in the limit, though this is not not rigorously shown here.

As a closely related issue, the Renyi divergence appears naturally in the 
study of normal approximation for densities $p_n$ of $Z_n$ in the form of 
non-uniform local limit theorems. Like in dimension one, denote by $\varphi$ 
the standard normal density in $\R^d$.

\vskip5mm
{\bf Theorem 1.3.} {\it Suppose that $D_\alpha(Z_n||Z)$ is finite for some 
$n=n_0$, and let the property $(1.4)$ be fulfilled. Then, for all $n$ 
large enough and for all $x \in \R^d$,
\be
|p_n(x) - \varphi(x)| \, \leq \, \frac{c}{\sqrt{n}}\, e^{-|x|^2/(2\alpha^*)}
\en
with some constant $c$ which does not depend on $n$. Moreover, the rate 
$1/\sqrt{n}$ on the right may be improved to $1/n$, provided that 
the distribution of $X$ is symmetric about the origin.
}

\vskip5mm
Thus, (1.5) is implied by the convergence $D_\alpha(Z_n||Z)\to 0$.
Non-uniform bounds in the normal approximation have been intensively studied in the
literature, cf. [Pe1-2], [I-L], [A1-2]. However, existing results start with weaker
hypotheses (e.g. moment assumptions) and either provide a polynomial error 
of approximation with respect to $x$ (such as $\frac{1}{1 + |x|^3}$), 
or deal with narrow zones contained in regions $|x| = o(\sqrt{n})$.

The paper consists of two parts.
In the first part results about the functional $D_\alpha$ are collected, 
including moment (exponential) inequalities and special properties of 
characteristic functions. Moreover, a number of remarkable algebraic properties 
of the $\chi^2$-distance will be derived.
They are related to the associated exponential series, the behavior under 
convolutions and heat semi-group transformations, and in higher dimensions -- 
to the super additivity of $\chi^2$ with respect to its marginals.
As a by-product, we establish the existence of densities in terms of the 
so-called normal moments. 
The second part is entirely devoted to the proof of Theorems 1.1-1.3. 
Employing an Edgeworth expansion for densities
(together with the results from the first part), this proof heavily relies 
on the tools of Complex Analysis. To simplify the presentation, almost 
all proofs will be stated for the one dimensional case, deffering 
the modifications needed to extend Theorems 1.1-1.3 to higher dimensions 
to separate sections. 

Thus the table of contents looks as follows:

\vskip5mm
\noindent
PART I: The $D_\alpha$ and $\chi^2$-divergence from the normal law

\vskip2mm
\noindent
2. Background on R\'enyi divergence

\noindent
3. Pearson-Vajda distances

\noindent
4. Basic exponential inequalities

\noindent
5. Laplace and Weierstrass transforms

\noindent
6. Connections with Fourier transform

\noindent
7. Exponential series

\noindent
8. Normal moments

\noindent
9. Behavior of R\'enyi divergence under convolutions

\noindent
10. Superadditivity of $\chi^2$ with respect to marginals

\vskip5mm
\noindent
PART II: The R\'enyi divergence in the central limit theorem

\vskip4mm
\noindent
11. Asymptotic expansions and lower bounds

\noindent
12. Necessity part in Theorem 1.2 ($d=1$)

\noindent
13. Pointwise upper bounds for convolutions of densities

\noindent
14. Sufficiency part in Theorem 1.2 ($d=1$)

\noindent
15. Non-uniform local limt theorem

\noindent
16. The multidimensional case

\noindent
17. Some examples and counter-examples

\noindent
18. Convolution of Bernoulli with Gaussian

\newpage
\vskip10mm
\noindent
{\large Part I: The $D_\alpha$ and $\chi^2$-divergence from the normal law}

\vskip10mm
\section{{\bf Background on R\'enyi Divergence}}
\setcounter{equation}{0}

\vskip2mm
\noindent
First let us briefly review some general properties of the R\'enyi divergences. More details can be 
found in the recent paper by van Erven and Harremo\"es [E-H]; cf. also [Le], [S], [G-S]. 

Let $(\Omega,\mu)$ be a measure space (with a $\sigma$-finite measure), and let 
$X$ and $Z$ be random elements with values in $\Omega$, having distributions $P$ and $Q$ 
with densities $p = \frac{dP}{d\mu}$, $q = \frac{dQ}{d\mu}$, respectively.
The following basic definitions go back to the work of R\'enyi [R]. 

\vskip5mm
{\bf Definition 2.1.} Let $0 < \alpha < \infty$, $\alpha \neq 1$.
The R\'enyi divergence of $P$ from $Q$ and the corresponding divergence power or relative 
Tsallis entropy of index $\alpha$ are the quantities
\bee
D_\alpha(X||Z)
 & = & 
D_\alpha(P||Q)  \, = \,
\frac{1}{\alpha - 1}\,\log \int\Big(\frac{p}{q}\Big)^\alpha q\, d\mu, \\
T_\alpha(X||Z) 
 & = &
T_\alpha(P||Q)  \, = \,
\frac{1}{\alpha - 1}\, \bigg[\int \Big(\frac{p}{q}\bigg)^\alpha q\, d\mu - 1\bigg].
\ene

The divergence $D_\alpha$ admits an axiomatic characterization via certain postulates. 
As a natural generalization of the Kullback-Leibler distance, the definition of $T_\alpha$ was 
introduced by Tsallis in [T] (within the so-called ``nonextensive thermostatistical formalism"), 
cf. also [B-T-P]. Both quantities are related by monotone transformations, namely
$$
D_\alpha = \frac{1}{\alpha - 1}\,\log\big(1 + (\alpha - 1)\,T_\alpha\big), \qquad
T_\alpha = \frac{1}{\alpha - 1}\, \big[e^{(\alpha - 1) D_\alpha} - 1\big].
$$
Thus, when they are small, these quantites are equivalent.
Both represent directional distances. In particular, $D_\alpha(P||Q) \geq 0$ and
$D_\alpha(P||Q) = 0$, if and only if $P=Q$.

The R\'enyi divergence with $0 < \alpha < 1$ posseses some unique features, like for example 
an obvious skew symmetry
$
D_\alpha(P||Q) = \frac{\alpha}{1-\alpha}\,D_\alpha(Q||P),
$
where the coefficient on the right is equal to 1 when 
$\alpha = \frac{1}{2}$. In this case, $D_\alpha$ represents 
a function of the square of the Hellinger metric:
$$
D_{1/2}(P||Q) = -2\log\Big(1 - \frac{1}{2}\,{\rm Hel}^2(P,Q)\Big).
$$
Another remarkable property is the equivalence of all $D_\alpha$ in this range: If
$0 < \alpha < \beta < 1$, then
$$
\frac{\alpha}{1 - \alpha}\,\frac{1-\beta}{1-\alpha}\, D_\beta(P||Q) \leq
D_\alpha(P||Q) \leq D_\beta(P||Q).
$$

When $\alpha \in (0,1)$ is fixed, $D_\alpha(P||Q)$ is a continuous function
of the tuple $(P,Q)$ with respect to the total variation distance in both coordinates. 
Conversely, it majorizes the total variation distance between $P$ and $Q$. 
Gilardoni [G] has shown that
$$
D_\alpha(P||Q) \geq \frac{\alpha}{2}\, \|P-Q\|_{\rm TV}^2.
$$
This extends the classical Pinsker inequality for the Kulback-Leibler distance
(when $\alpha=1$), with best constant due to Csisz\'ar, cf.  [Pi], [Cs].

The following general property is important for comparing the R\'enyi
divergence with different values of $\alpha$.

\vskip5mm
{\bf Proposition 2.2.} {\it For all probability measures $P$ and $Q$ on 
$\Omega$, the functions $\alpha \rightarrow D_\alpha(P||Q)$ and 
$\alpha \rightarrow T_\alpha(P||Q)$ are non-decreasing.
}

\vskip5mm
The monotonicity of $D_\alpha$ is discussed in [E-H], Theorem 3.
As for $T_\alpha$, let $0 < \alpha < \beta$, $\alpha, \beta \neq 1$.
The functions $c \rightarrow e^{ct_0} - 1$ with fixed $t_0 \geq 0$ and 
$t \rightarrow \frac{e^{ct} - 1}{t}$ are non-decreasing in $c \geq 0$ and $t>1$, respectively. 
Hence, in case $\alpha > 1$, we get, using monotonicity of $D_\alpha$,
\bee
T_\alpha(P||Q) 
 & = &
\frac{1}{\alpha - 1}\, \big[e^{(\alpha - 1) D_\alpha(P||Q)} - 1\big]
 \ \leq \ 
\frac{1}{\alpha - 1}\, \big[e^{(\alpha - 1) D_\beta(P||Q)} - 1\big] \\
 & \leq &
\frac{1}{\beta - 1}\, \big[e^{(\beta - 1) D_\beta(P||Q)} - 1\big]
 \ = \ T_\beta(P||Q).
\ene
In case $\alpha < 1$, we use the property that the function
$c \rightarrow 1-e^{-ct_0}$ is non-decreasing in $c \geq 0$, while
$t \rightarrow \frac{1-e^{-ct}}{t}$ is non-increasing on the half-axis 
$-\infty < t < 1$. This yields
\bee
T_\alpha(P||Q) 
 & = &
\frac{1}{1-\alpha}\, \big[1-e^{-(1-\alpha) D_\alpha(P||Q)}\big]
 \ \leq \
\frac{1}{1-\alpha}\, \big[1-e^{-(1-\alpha) D_\beta(P||Q)}\big] \\
 & \leq &
\frac{1}{1-\beta}\, \big[1-e^{-(1-\beta) D_\beta(P||Q)}\big]
 \ = \ T_\beta(P||Q).
\ene

The values $0<\alpha<1$ and $1 < \alpha < \infty$, for which the R\'enyi divergence 
was defined explicitly, are called simple. The monotonicity of $D_\alpha(P||Q)$ with 
respect $\alpha$ allows to extend this function to the missing values $\alpha = 0$, 
$\alpha = 1$ and $\alpha = \infty$, which are called extended values:
\bee
D_0(P||Q) & = & \lim_{\alpha \downarrow 0} D_\alpha(P||Q), \qquad
D_\infty(P||Q) \ = \ \lim_{\alpha \rightarrow \infty} D_\alpha(P||Q), \\
D_1(P||Q) & = & \lim_{\alpha \uparrow 1} D_\alpha(P||Q).
\ene
It is easy to check that
$D_0(P||Q) = -\log\, Q\{p(x) > 0\}$ and
$D_\infty(P||Q) \, = \, \log\, {\rm ess\,sup}_P\, \frac{p(x)}{q(x)}$
with the convention that $0/0 = 0$.

The extended index $\alpha=0$ may be used to characterize an absolute 
continuity or singularity of two given probability distributions:
$D_0(P||Q) = 0$, if and only if $Q$ is absolutely continuous with respect to $P$,
and $D_0(P||Q) = \infty$, if and only if $P$ and $Q$ are orthogonal to each other.
This can be illustrated by the Gaussian dichotomy -- the property saying that any 
two Gaussian measures are either absolutely continuous to each other or orthogonal, 
cf. [S], p. 366.

The extended index $\alpha=1$ leads to the Kullback-Leibler distance
$$
D(X||Z) = D(P||Q) = \int p\,\log \frac{p}{q}\ d\mu,
$$
also known as the relative entropy or the informational divergence. 
Motivated by works of Shannon and Wiener on communication engineering, 
this quantity was introduced by Kullback and Leibler [K-L] under the name 
``the information of $P$ relative to $Q$'' (though using a different notation). 
Note that in this case $D_1 = T_1 = D$. It follows immediately that $D(P||Q) = \infty$, 
if $P$ is not absolutely continuous with respect to $Q$.

As was already mentioned, in the particular case $\alpha = 2$, we arrive at the definition 
of the quadratic Renyi divergence and the quadratic Renyi divergence power also known
as the $\chi^2$-distance:
$$
D_2(X||Z) \, = \log \int \frac{p^2}{q}\, d\mu, \qquad
\chi^2(X,Z) \, = \, T_2(X||Z) \, = \int \frac{p^2}{q}\, d\mu - 1.
$$
In all cases, by the Csisz\'ar-Pinsker inequality for $\alpha=1$, we have the relations
$$
\frac{1}{2}\,\|P-Q\|_{\rm TV}^2 \leq D(X||Z) \leq D_2(X||Z) \leq \chi^2(X,Z).
$$

Another important property of these distances is the contractivity under mappings.

\vskip5mm
{\bf Proposition 2.3.} {\it For any measurable map $S$ from $\Omega$ 
to any measurable space $\Omega'$,
\be
D_\alpha(S(X)||S(Z)) \leq D_\alpha(X||Z) \qquad (\alpha \geq 1).
\en
}

\vskip2mm
{\bf Proof.} Suppose that $D_\alpha(X||Z)$ is finite, so that the distribution $P$ is 
absolutely continuous with respect $Q$. Introducing $\xi = p/q$, $\beta = \alpha/(\alpha - 1)$
with $\alpha>1$, one may write
\bee
\big(1 + (\alpha - 1)\, T_\alpha(X||Z)\big)^{1/\alpha}
 & = &
(\E_Q\, \xi^\alpha)^{1/\alpha} \ = \ \sup_{\E_Q \eta^\beta \leq 1} \E_Q\,\xi \eta \\
 & = & 
\sup_{\E_Q \eta^\beta \leq 1} \E_P\,\eta \ = \ 
\sup_{\E\, \eta(Z)^\beta \leq 1} \E\,\eta(X),
\ene
that is,
\be
1 + (\alpha - 1)\, T_\alpha(X||Z) \, = \, \sup_{\E\, \eta(Z)^\beta \leq 1} 
\big(\E\,\eta(X)\big)^\alpha,
\en
where the sup is taken over all measurable functions $\eta:\Omega \rightarrow \R_+$
such that $\E\, \eta(Z)^\beta \leq 1$. Similarly
$$
1 + (\alpha - 1)\, T_\alpha(S(X),S(Z)) \, = \, 
\sup_{\E\, \eta(S(Z))^\beta \leq 1} \big(\E\,\eta(S(X))\big)^\alpha
 \, = \, \sup_{\E\, \eta'(Z)^\beta \leq 1} (\E\,\eta'(X))^\alpha
$$
where the second supremum on the right has been restricted to the class of functions of 
the form $\eta' = \eta(S)$. Hence, this supremum does not exceed the right-hand side of (2.2), 
thus proving (2.1) for $T_\alpha$.
\qed

\vskip5mm
The property (2.1) is closely related to the so called data processing inequality in 
Information Theory, namely
$$
D_\alpha(P_{{\mathfrak A}} || Q_{{\mathfrak A}}) \leq D_\alpha(P||Q),
$$
where $P_{{\mathfrak A}}$ and $Q_{{\mathfrak A}}$ denote restrictions of the measures
$P$ and $Q$ to an arbitrary $\sigma$-subalgebra $\mathfrak A$ in $\Omega$
(cf. [E-H], Theorem 1).

\vskip10mm
\section{{\bf Pearson-Vajda Distances}}
\setcounter{equation}{0}

\vskip2mm
\noindent
Writing
$
\chi^2(X,Z) \,= \, \int \frac{|p-q|^2}{q}\, d\mu,
$
the $\chi^2$-distance may be regarded as a particular member in the family of
Pearson-Vajda distances [N], descibed below.

\vskip5mm
{\bf Definition 3.1.} For $\alpha \geq 1$, the $\chi_\alpha$-distance of $P$ from $Q$
is defined by
$$
\chi_\alpha(X,Z) = \chi_\alpha(P,Q) = 
\int \Big|\frac{p-q}{q}\Big|^\alpha q\, d\mu = 
\big\|p-q\big\|_{L^\alpha(\frac{1}{q^{\alpha - 1}}\,d\mu)}^\alpha.
$$

\vskip5mm
As in the previous section, here $X$ and $Z$ denote random elements in $(\Omega,\mu)$, 
having distributions $P$ and $Q$  with densities $p = \frac{dP}{d\mu}$, $q = \frac{dQ}{d\mu}$.
The quantity $\chi_\alpha(P,Q)$ (which is often denoted $\chi^\alpha$)
does not depend on the choice of the dominating measure $\mu$.

Clearly, the function $\chi_\alpha^{1/\alpha}$ is non-decreasing in $\alpha$, and when 
$\alpha = 1$, we arrive at the total variation distance between $P$ and $Q$.

For our further purpose, it will be useful to relate the R\'enyi divergence power $T_\alpha$
to $\chi_\alpha$. Both quantities are metrically equivalent, as seen by the following
elementary observation.

\vskip5mm
{\bf Proposition 3.2.} {\it For all $\alpha > 1$, 
\be
T_\alpha \, \leq \, \frac{1}{\alpha - 1}\, 
\Big[\big(1 + \chi_\alpha^{1/\alpha}\big)^\alpha - 1\Big],
\en
where $T_\alpha = T_\alpha(P||Q)$ and $\chi_\alpha = \chi_\alpha(P||Q)$.
Conversely,
\be
T_\alpha \, \geq \, \frac{3}{16}\,\min\{\chi_\alpha,\chi_\alpha^{2/\alpha}\} \ \ 
(1 < \alpha \leq 2), \ \qquad \ T_\alpha\, \geq \, 3^{-\alpha}\, \chi_\alpha \ \ (\alpha \geq 2).
\en
}

\vskip2mm
{\bf Proof.}
By the triangle inequality in $L^\alpha(\frac{1}{q^{\alpha - 1}}\,d\mu)$,
\bee
\chi_\alpha^{1/\alpha} \ = \ 
 \|p-q\|_{L^\alpha(\frac{1}{q^{\alpha - 1}}\,d\mu)} 
 & \geq &
\Big|\,\|p\|_{L^\alpha(\frac{1}{q^{\alpha - 1}}\,d\mu)} -
\|q\|_{L^\alpha(\frac{1}{q^{\alpha - 1}}\,d\mu)}\Big| \\
 & = &
\bigg(\int \Big(\frac{p}{q}\Big)^\alpha q\, d\mu\bigg)^{1/\alpha} - 1
 \, = \,
\big(1 + (\alpha - 1)\,T_\alpha\big)^{1/\alpha} - 1,
\ene
which proves (3.1).

To argue in the opposite direction, put $\xi = p/q$. Since $dQ = q\,d\mu$, we may write
$$
T_\alpha = \frac{1}{\alpha - 1}\, \big[\E\,\xi^\alpha  - 1\big], \qquad
\chi_\alpha = \E\,|\xi - 1|^\alpha,
$$
where the expectations are taken on the probability space $(\Omega,Q)$. We have
$\xi \geq 0$ and $\E \xi = 1$. Consider the random variable $\eta = \xi - 1 \geq -1$
and the function
$$
\psi(t) = \E\,(1 + t\eta)^\alpha  - 1, \qquad t \geq 0,
$$
so that $\psi(1) = \E\,\xi^\alpha  - 1$. This function is differentiable in $t>0$, with
continuous derivatives
$$
\psi'(t) = \alpha\, \E\,\eta (1 + t\eta)^{\alpha - 1} , \qquad 
\psi''(t) = \alpha (\alpha - 1)\, \E\,\eta^2  (1 + t\eta)^{\alpha - 2}.
$$
Since $\psi(0) = \psi'(0) = 0$, by the Taylor integral formula,
$$
\psi(1) = \int_0^1 (1-t) \psi''(t)\,dt = 
\alpha (\alpha - 1)\, \E\,\eta^2 \int_0^1 (1-t) (1 + t\eta)^{\alpha - 2}\,dt.
$$

{\it Case} $1 < \alpha \leq 2$. Since the function $t \rightarrow (1 + t\eta)^{\alpha - 2}$
is convex on $(0,\infty)$, Jensen's inequality with respect to the probability measure
$d\nu(t) = 2(1-t)\,dt$ on $(0,1)$ yields
\bee
\int_0^1 (1-t) (1 + t\eta)^{\alpha - 2}\,dt 
 & = &
\frac{1}{2}\,\int(1 + t\eta)^{\alpha - 2}\,d\nu(t) \\
 & \geq &
\frac{1}{2}\,\bigg(1 + \eta \int t\,d\nu(t)\bigg)^{\alpha - 2} \ = \
\frac{1}{2}\,\Big(1 + \frac{1}{3}\,\eta\Big)^{\alpha - 2}.
\ene
Therefore,
$$
\psi(1) \geq \frac{1}{2}\, \alpha (\alpha - 1)\, 
\E\,\eta^2 \Big(1 + \frac{1}{3}\,\eta\Big)^{\alpha - 2}.
$$

On the set $A = \{|\eta| \leq 1\}$, the expression $\eta^2 (1 + \frac{1}{3}\,\eta)^{\alpha - 2}$
is bounded from below by $(\frac{3}{4})^{2 - \alpha}\, \eta^2$, and on the set 
$B = \{\eta > 1\}$ by
$\eta^2 \cdot (\frac{4}{3}\,\eta)^{\alpha - 2} = (\frac{3}{4})^{2 - \alpha}\, \eta^\alpha$. 
Hence
$$
\psi(1) \geq \frac{1}{2}\, \alpha (\alpha - 1)\,  \Big(\frac{3}{4}\Big)^{2 - \alpha}\,
\E\,\big(\eta^2\,1_A + \eta^\alpha 1_B\big).
$$
For our range of $\alpha$'s we may simply use $\alpha \,(\frac{3}{4})^{2 - \alpha} \geq \frac{3}{4}$, 
so that
$$
\psi(1) \geq \frac{3}{8}\,(\alpha - 1)\,
\E\,\big(\eta^2\,1_A + \eta^\alpha 1_B\big).
$$
By Markov's inequality,
$$
\frac{1}{\P(A)}\,\E\,\eta^2\,1_A \geq 
\bigg(\frac{1}{\P(A)}\,\E\,|\eta|^\alpha\,1_A\bigg)^{2/\alpha}, 
$$
so
$\E\,\eta^2\,1_A \geq (\E\,|\eta|^\alpha\,1_A)^{2/\alpha}$ and thus
$$
\E\,\big(\eta^2\,1_A + \eta^\alpha\, 1_B\big) \geq U = u_0^{2/\alpha} + u_1, \quad {\rm where} \ \
u_0 = \E\,|\eta|^\alpha\,1_A, \  \ u_1 = \E\,|\eta|^\alpha\,1_B.
$$
Fixing the value $u = u_0 + u_1$, in case $u_1 \geq \frac{1}{2}\,u$ we have
$U \geq u_1 \geq \frac{1}{2}\,u$, while in case $u_0 \geq \frac{1}{2}\,u$ we have
$U \geq u_0 \geq (\frac{1}{2}\,u)^{2/\alpha} \geq \frac{1}{2}\,u^{2/\alpha}$. In both cases,
$U \geq \frac{1}{2}\,\min(u,u^{2/\alpha})$, that is,
$$
\E\,\big(\eta^2\,1_A + \eta^\alpha\, 1_B\big) \geq 
\frac{1}{2}\, \min\{\E\,|\eta|^\alpha, (\E\,|\eta|^\alpha)^{2/\alpha}\}.
$$
As a result,
$$
T_\alpha \ = \ \frac{1}{\alpha - 1}\,\psi(1) 
 \, \geq \,
\frac{3}{16}\, \min\{\E\,|\eta|^\alpha, (\E\,|\eta|^\alpha)^{2/\alpha}\}
 \, = \,
\frac{3}{16}\,\min\{\chi_\alpha,\chi_\alpha^{2/\alpha}\},
$$
which yields the first inequality in (3.2).

{\it Case} $\alpha > 2$. Let us return to the Taylor integral formula
$$
\psi(1) = \alpha (\alpha - 1)\, \E\,\eta^2 \int_0^1 (1-t) (1 + t\eta)^{\alpha - 2}\,dt,
$$
where we now restrict integration to the interval $(\frac{1}{3},\frac{2}{3})$ to get
$$
\psi(1) \geq \frac{\alpha (\alpha - 1)}{3}\, \E\,\eta^2 \int_{1/3}^{2/3} (1 + t\eta)^{\alpha - 2}\,dt.
$$
Since $\eta \geq - 1$, in case $\eta \leq 0$, we have
$1 + t\eta \geq 1 + \frac{2}{3}\,\eta \geq -\frac{1}{3}\,\eta$. 
In case $\eta \geq 0$, we similarly have $1 + t\eta \geq t \eta \geq \frac{1}{3}\,\eta$. In both cases,
$1 + t\eta \geq \frac{1}{3}\,|\eta|$, hence
$$
\psi(1) \geq \frac{\alpha (\alpha - 1)}{3}\, 
\E\,\eta^2 \int_{1/3}^{2/3} \Big(\frac{1}{3}\,|\eta|\Big)^{\alpha - 2}\,dt = 
\frac{\alpha (\alpha - 1)}{3^\alpha}\, \E\,|\eta|^\alpha,
$$
and therefore
$$
T_\alpha \, = \, \frac{1}{\alpha - 1}\,\psi(1) \, \geq \,
\frac{\alpha - 1}{3^\alpha}\, \E\,|\eta|^\alpha \, = \, \frac{\alpha - 1}{3^\alpha}\, \chi_\alpha,
$$
\qed

\vskip10mm
\section{{\bf Basic Exponential Inequalities}}
\setcounter{equation}{0}

\vskip2mm
\noindent
We now focus on the particular case, where $\Omega = \R$ is the real line 
with Lebesgue measure $\mu$, and where $Z \sim N(0,1)$ is a standard normal random 
variable, i.e., with density 
$$
\varphi(x) = \frac{1}{\sqrt{2\pi}}\,e^{-x^2/2}, \qquad x \in \R.
$$
Given a random variable $X$,  the R\'enyi divergence and the Tsallis distance of index $\alpha>1$
with respect to $Z$ are then given by the formulas
$$
(\alpha - 1)\,D_\alpha(X||Z) = 
\log \int_{-\infty}^\infty \frac{p(x)^\alpha}{\varphi(x)^{\alpha - 1}}\, dx,  \qquad
(\alpha - 1)\,T_\alpha(X||Z) = 
\int_{-\infty}^\infty \frac{p(x)^\alpha}{\varphi(x)^{\alpha - 1}}\, dx - 1,
$$
where $p$ is density of $X$. 
If the distribution of $X$ is not absolutely continuous
with respect to $\mu$, then we automatically have $D_\alpha(X||Z) = T_\alpha(X||Z) = \infty$.
These quantities are finite, if, for example, $p$ is bounded and $\E\,e^{(\alpha - 1) X^2/2} < \infty$.
In fact, the finiteness of $D_\alpha(X||Z)$ or $T_\alpha(X||Z)$ implies a similar property.
In the sequel, we put
$$
\beta = \alpha^* = \frac{\alpha}{\alpha - 1}.
$$

\vskip5mm
{\bf Proposition 4.1.} {\it If $T_\alpha =T_\alpha(X||Z) < \infty$, then 
$X$ must have an absolutely continuous distribution with 
$$
\E\,e^{cX^2} \, \leq \, \frac{C}{(1-2\beta c)^{\frac{1}{2\beta}}} \quad {for \ all} \ \ 
c < \frac{1}{2\beta},
$$
where $C = \big(1 + (\alpha - 1)T_\alpha\big)^{1/\alpha}$.
It is possible that $T_\alpha < \infty$, while $\E\,e^{\frac{1}{2\beta} X^2} =\infty$.
}

\vskip5mm
In particular, if $T_\alpha$ is finite, $X$ must finite moments of any order.

\vskip5mm
{\bf Proof.} Let $X$ have density $p$ such that the integral
$C = \big(\int_{-\infty}^\infty \frac{p(x)^\alpha}{\varphi(x)^{\alpha - 1}}\, dx\big)^{1/\alpha}$ 
is finite. By the H\"older inequality with dual exponents $(\beta,\alpha)$,
\bee
\E\,e^{cX^2}
 & = &
\int_{-\infty}^\infty
\frac{p(x)}{\varphi(x)^{1/\beta}}\cdot e^{cx^2}\varphi(x)^{1/\beta}\, dx \\
 & \leq &
C\, \bigg(\int_{-\infty}^\infty e^{\beta cx^2}\,\varphi(x)\, dx\bigg)^{1/\beta} 
\ = \ \frac{C}{(1-2\beta c)^{\frac{1}{2\beta}}}\, .
\ene
This proves the first assertion. For the second assertion, one may consider a density of the form
$
p(x) = \frac{a}{1 + |x|}\,e^{-\frac{1}{2\beta} x^2},
$
where $a$ is a normalizing constant. Then
$T_\alpha < \infty$ and $\E\,e^{\frac{1}{2\beta} X^2} = \infty$. 
\qed

\vskip5mm
As an alternative (although almost equivalent) variant of Proposition 4.1, we also have:

\vskip5mm
{\bf Proposition 4.2.} {\it If $T_\alpha =T_\alpha(X||Z) < \infty$, then for all $t \in \R$,
\be
\E\,e^{tX} \leq Ce^{\beta t^2/2},
\en
where $C = \big(1 + (\alpha - 1)T_\alpha\big)^{1/\alpha}$.
In particular, for any $r \geq 0$, 
$$
\P\{X \geq r\} \leq C e^{-\frac{1}{2\beta} r^2}.
$$
}

\vskip2mm
Indeed, arguing as before, if $p$ is density of $X$,
\bee
\E\,e^{tX}
 & = &
\int_{-\infty}^\infty p(x)\,e^{tx}\, dx \\
 & = &
\int_{-\infty}^\infty
\frac{p(x)}{\varphi(x)^{1/\beta}}\cdot e^{tx}\varphi(x)^{1/\beta}\, dx
 \, \leq \, 
C\, \bigg(\int_{-\infty}^\infty e^{\beta tx}\,\varphi(x)\, dx\bigg)^{1/\beta} 
= \, Ce^{\beta t^2/2}.
\ene

This bound cannot be deduced from the bound of Proposition 4.1. In fact, the coefficient 
$C$ in (4.1) may be chosen to be smaller than 1 for large values of $|t|$. 
The next assertion will be one of the steps needed in the proof of the sufficiency part of 
Theorems 1.1-1.2.

\vskip5mm
{\bf Proposition 4.3.} {\it If $T_\alpha(X||Z) < \infty$, then
$$
\lim_{|t| \rightarrow \infty} \E\,e^{tX}\,e^{-\beta t^2/2} = 0.
$$
}

\vskip2mm
{\bf Proof.} Let $p$ be the density of $X$ and write
$
\E\,e^{tX} = \int_{-\infty}^{\infty} e^{tx}\,p(x)\,dx.
$
Here integration over the positive half-axis may be splitted into the two intervals.
First, given $t>0$, by the H\"older inequality,
\bee
\int_0^{\beta t/2} e^{tx}\,p(x)\,dx
 & = &
\int_0^{\beta t/2} p(x)\, e^{\frac{x^2}{2\beta}}\cdot e^{tx - \frac{x^2}{2\beta}}\,dx \\ 
 & \le &
\bigg(\int_{-\infty}^\infty p(x)^\alpha\,e^{\frac{\alpha x^2}{2\beta}}\,dx\bigg)^{1/\alpha}
\bigg(\int_0^{\beta t/2} e^{\beta tx - \frac{x^2}{2}}\,dx\bigg)^{1/\beta}\\
 & \le &
\frac{1}{(2\pi)^{1/(2\beta)}}\, \big(1+(\alpha - 1)T_\alpha(X||Z)\big)^{1/\alpha}\, 
\Big(\frac{\beta t}{2}\Big)^{1/\beta}\,e^{3\beta t^2/8},
\ene
where we used the monotonicity of $\beta tx - \frac{1}{2}\,x^2$ in the interval $0 \leq x \leq \beta t$
(in order to estimate the last integral). 
Similarly,
$$
\int_{\beta t/2}^\infty p(x)\, e^{\frac{x^2}{2\beta}} \cdot e^{tx - \frac{x^2}{2\beta}}\,dx \,\le \,
\bigg(\int_{\beta t/2}^\infty p(x)^\alpha\,e^{\frac{\alpha x^2}{2\beta}}\,dx\bigg)^{1/\alpha}
\bigg(\int_{-\infty}^\infty e^{\beta tx - \frac{x^2}{2}}\,dx\bigg)^{1/\beta} \, \leq \,
\delta(t)\, e^{\beta t^2/2}
$$
with $\delta(t)\to 0$ as $t\to\infty$. Collecting these bounds, we get
$$
\E\, e^{tX}\,1_{\{X>0\}}\, e^{-\beta t^2/2} \, \le \,
(2\pi)^{-1/(2\beta)}\,\big(1+(\alpha - 1)T_\alpha(X||Z)\big)^{1/\alpha}\, 
\Big(\frac{\beta t}{2}\Big)^{1/\beta}\, e^{-\beta t^2/8} + 
\delta(t) \to 0.
$$
Since also $\E\, e^{tX}\,1_{\{X<0\}} \rightarrow 0$ as $t \rightarrow \infty$, the
conclusion follows.
\qed

\vskip10mm
\section{{\bf Laplace and Weierstrass Transforms}}
\setcounter{equation}{0}

\vskip2mm
\noindent
Although in general the critical constant in the exponent $c = 1/(2\beta)$ cannot be included
in the statement of Proposition 4.1, this turns out possible for suffiently many normalized convolutions 
of the distribution of $X$ with itself. Given  independent copies $X_1,\dots,X_n$
of $X$, here we consider ``Gaussian" moments for the normalized sums
$$
Z_n = \frac{X_1 + \dots + X_n}{\sqrt{n}}.
$$
The following statement is crucial both in the necessity and sufficiency parts of the proof
of Theorems 1.1-1.2. We always assume that $Z \sim N(0,1)$.

\vskip5mm
{\bf Proposition 5.1.} {\it If $T_\alpha = T_\alpha(X||Z) < \infty$, then
$\E\,e^{\frac{1}{2\beta} Z_n^2} < \infty$ for all $n \geq \alpha$, and
\be
\E\,e^{\frac{1}{2\beta} Z_n^2} \, \leq \, 3^n \big(1 + (\alpha - 1)T_\alpha\big)^{\frac{n}{\alpha}}.
\en
Moreover, putting $\chi_\alpha = \chi_\alpha(X,Z)$, we have
\be
\big|\E\,e^{\frac{1}{2\beta} Z_n^2} - \E\,e^{\frac{1}{2\beta} Z^2}\big|  \, \leq \,
3^n \Big(\big(1 + \chi_\alpha^{1/\alpha}\big)^n - 1\Big).
\en
}

\vskip2mm
Thus, when $X$ is close to $Z$ in the sense of the Pearson-Vajda distance, we also obtain 
closeness of the corresponding Gaussian moments of $Z_n$ and $Z$ with fixed $n \geq \alpha$. 
Recall that $\chi_\alpha$ in (5.2) can be estimated from above in terms of $T_\alpha$ 
according to Proposition 3.2 (while these distances coincide in case $\alpha = 2$).

As for the inequality (5.1), one may equivalently rephrase it in terms of the Laplace transform 
of the distribution of $Z_n$. Let us state one immediate corrollary.

\vskip5mm
{\bf Corollary 5.2.} {\it  Let $T_\alpha =T_\alpha(X||Z)$ be finite. Then the function 
$\psi(t) = \E\,e^{t X}\,e^{-\beta t^2/2}$ is integrable with any power $n \geq \alpha$, 
and up to some $n$-dependent constant $c_n$,
\be
\int_{-\infty}^{\infty} \psi(t)^n\,dt \, \le \, c_n \big(1 + (\alpha - 1)T_\alpha\big)^{\frac{n}{\alpha}}.
\en
}

The argument uses the contractivy properties of the Weierstrass transform, which is defined by 
the equality
$$
W_t u(x) = \frac{1}{\sqrt{2\pi t}} \int_{-\infty}^\infty e^{-\frac{(x-y)^2}{2t}}\,u(y)\,dy, \qquad
x \in \R, \ t > 0.
$$

For short in the sequel we denote by $L^\alpha$ the Lebesgue space $L^\alpha(\R,dx)$ of all 
measurable functions on the real line with finite norm
$$
\|u\|_\alpha = \Big(\int_{-\infty}^\infty |u(x)|^\alpha\,dx\Big)^{1/\alpha}, \qquad \alpha \geq 1,
$$
with usual convention $\|u\|_\infty = {\rm ess\,sup}_x\, |u(x)|$.

We refer an interested reader to [H-W] for a detail account on the Weierstrass 
transform, and here only mention one property.
Since $W_t u$ represents the convolution of $u$, namely, with the Gaussian density
$\varphi_t(x) =  \frac{1}{\sqrt{2\pi t}}\, e^{-x^2/(2t)}$, we have, by Jensen's
inequality, $\|W_t u\|_\alpha \leq \|u\|_\alpha$ for all $\alpha \geq 1$ and $t > 0$. That is, $W_t$
acts as a contraction from $L^\alpha$ to $L^\alpha$. 

This implies that $W_t$ is a bounded operator from $L^\alpha$ to $L^\gamma$ with any 
$\gamma > \alpha$. Indeed, by H\"older's inequality, 
$|W_t u(x)| \leq \|\varphi_t\|_\beta\, \|u\|_\alpha$ $(\beta = \alpha^*)$ , and since
$
\|\varphi_t\|_\beta = (2\pi t)^{-1/(2\alpha)}\,\beta^{-1/(2\beta)},
$
we get
$$
\|W_t u\|_\infty \, \leq \, (2\pi t)^{-1/(2\alpha)}\,\beta^{-1/(2\beta)}\,\|u\|_\alpha.
$$
More generally, given $\alpha < \gamma < \infty$, we have
$$
\int_{-\infty}^\infty |W_t u(x)|^\gamma\,dx
 \, = \,
\int_{-\infty}^\infty |W_t u(x)|^{\gamma - \alpha}\,|W_t u(x)|^\alpha\,dx
 \, \leq \,
(2\pi t)^{\frac{\alpha - \gamma}{2\alpha}}\,\beta^{\frac{\alpha - \gamma}{2\beta}}\,
\|u\|_\alpha^\gamma.
$$
Hence 
\be
\|W_t u\|_\gamma \, \leq \,
(2\pi t)^{\frac{\alpha - \gamma}{2\gamma\alpha}}\,
\beta^{\frac{\alpha - \gamma}{2 \gamma\beta}}\,\|u\|_\alpha, \qquad
\alpha \leq \gamma \leq \infty.
\en
In fact, since 
$\frac{\alpha - \gamma}{\gamma\alpha} = \frac{1}{\gamma} - \frac{1}{\alpha}$ 
may vary from zero to $-\frac{1}{\alpha}$, 
the latter bound can be made independent of $\gamma$, namely, in the indicated range
$$
\|W_t u\|_\gamma \, \leq \, \max\big\{1, (2\pi t)^{-\frac{1}{2\alpha}}\big\}\,\|u\|_\alpha.
$$

The inequality (5.4) is what we need for the proof of Proposition 5.1.

\vskip5mm
{\bf Proof of Proposition 5.1.} Let $p$ be the density of $X$. The Weierstrass transform 
can be applied to the function
$$
u(x) = \varphi(x)^{-1/\beta}\, p(x),
$$
which has finite norm $\|u\|_\alpha = (1 + (\alpha - 1)T_\alpha)^{1/\alpha}$. Putting
$\bar x = \frac{1}{n}\,(x_1 + \dots + x_n)$,
the expectation we have to estimate is
\bee
\E\,e^{\frac{1}{2\beta} Z_n^2}
 & = &
\int_{\R^n} e^{\frac{n}{2\beta}\,\bar x^2}\,p(x_1)\dots p(x_n)\,dx_1 \dots dx_n \\
 & = &
(2\pi)^{-\frac{n-1}{2\beta}} \int_{\R^n} 
\exp\Big\{\frac{n}{2\beta}\,\bar x^2 - \frac{1}{2\beta}\,(x_1^2 + \dots + x_n^2)\Big\}\, 
u(x_1)\dots u(x_n)\,dx_1 \dots dx_n \\
 & = &
(2\pi)^{-\frac{n-1}{2\beta}} \int_{\R^n} 
\exp\Big\{-\frac{1}{4\beta n}\,\sum_{i=1}^n Q_i\Big\}\, 
u(x_1)\dots u(x_n)\,dx_1 \dots dx_n,
\ene
where $Q_i = \sum_{j = 1}^n (x_i - x_j)^2$. First,
we apply H\"older's inequality and put $t = 2\beta$, to get
\bee
\E\,e^{\frac{1}{2\beta} Z_n^2} 
 & \leq &
(2\pi)^{-\frac{n-1}{2\beta}} \prod_{i=1}^n \bigg(\int_{\R^n} 
\exp\Big\{-\frac{1}{4\beta}\,\sum_{i=1}^n Q_i\Big\}\, 
u(x_1)\dots u(x_n)\,dx_1 \dots dx_n\bigg)^{1/n} \\ 
 & = &
(2\pi)^{-\frac{n-1}{2\beta}} \, (2\pi t)^{\frac{n-1}{2}} 
\int_{-\infty}^\infty (W_t u(x_1))^{n-1}\, u(x_1)\, dx_1,
\ene
where on the second step, inside the $i$-th integral in the product we performed 
the integration over the variables $x_j$, $j \neq i$, which yielded the value 
$(2\pi t)^{\frac{n-1}{2}}\,(W_t u(x_i))^{n-1}$. By H\"older's inequality once more, 
and applying (5.4) with $\gamma = \beta(n-1)$, which satisfies $\gamma \geq \alpha$ 
due to the assumption $n \geq \alpha$, we see that the last one dimensional integral 
does not exceed
\bee
\Big(\int_{-\infty}^\infty (W_t u(x_1))^\gamma\,dx_1\Big)^{\frac{1}{\beta}}\, 
\|u\|_\alpha
 & = &
\|W_t u\|_\gamma^{n-1}\, \|u\|_\alpha \\
 & \leq &
\big((2\pi t)^{\frac{\alpha - \gamma}{2\gamma\alpha}}\,
\beta^{\frac{\alpha - \gamma}{2 \gamma\beta}}\, \|u\|_\alpha\big)^{n-1}\, \|u\|_\alpha.
\ene
Hence
$
\E\,e^{\frac{1}{2\beta} Z_n^2} \, \leq \, c_{n,\alpha} \|u\|_\alpha^n
$
with constant
\bee
c_{n,\alpha}
 & = & 
(2\pi)^{-\frac{n-1}{2\beta}}\, (2\pi t)^{\frac{n-1}{2}}\,
(2\pi t)^{\frac{n-1}{2}\,\frac{\alpha - \gamma}{\gamma \alpha}}\, 
\beta^{\frac{n-1}{2}\,\frac{\alpha - \gamma}{\gamma\beta}}\\
 & = &
(2\pi)^{\frac{1}{2\beta}}\,t^{\frac{n}{2\beta}\,}\beta^{\frac{\alpha - n}{2\beta}}
 \ = \ 
(2\pi)^{\frac{1}{2\beta}}\, 2^{\frac{n}{2\beta}}\,\beta^{\frac{1}{2(\beta - 1)}} \, < \, 3^n.
\ene
This proves (5.1). 
It is also interesting to note that $c_{n,\alpha} \rightarrow 1$ as $\alpha \rightarrow 1$.

Obviously, this argument can easily be extended to not necessarily equal positive functions.
Namely, for the integral
$$
I = I(p_1,\dots,p_n) = \int_{\R^n} 
e^{\frac{n}{2\beta}\,\bar x^2}\,p_1(x_1)\dots p_n(x_n)\,dx_1 \dots dx_n
$$
we similarly obtain
\bee
 |I|
 & \leq &
(2\pi)^{-\frac{n-1}{2\beta}} \prod_{i=1}^n \bigg(\int_{\R^n} 
\exp\Big\{-\frac{1}{4\beta}\,\sum_{i=1}^n Q_i\Big\}\, 
|u(x_1)|\dots |u(x_n)|\,dx_1 \dots dx_n\bigg)^{1/n} \\ 
 & = &
(2\pi)^{-\frac{n-1}{2\beta}} \, (2\pi t)^{\frac{n-1}{2}} 
\prod_{i=1}^n \bigg( \int_{-\infty}^\infty |u_i(x_i)|\, \prod_{j \neq i}\,
(W_t |u_j|)(x_i)\, dx_i\bigg)^{1/n},
\ene
where $u_j = \varphi^{- 1/\beta}\, p_j$. An application of H\"older's inequality together with (5.4)
allows one to estimate the last integral by
\bee
\|u_i\|_\alpha \prod_{j \neq i} \|W_t |u_j|\|_\gamma 
 & \leq & 
\|u_i\|_\alpha \prod_{j \neq i} \,
(2\pi t)^{\frac{\alpha - \gamma}{2\gamma \alpha}}\, 
\beta^{\frac{\alpha - \gamma}{2\gamma\beta}}\,\|u_j\|_\alpha \\
 & = &
(2\pi t)^{\frac{n-1}{2}\,\frac{\alpha - \gamma}{\gamma \alpha}}\, 
\beta^{\frac{n-1}{2}\,\frac{\alpha - \gamma}{\gamma\beta}}\,
\|u_1\|_\alpha \dots \|u_n\|_\alpha.
\ene
This leads to
\be
|I(p_1,\dots,p_n)| \leq c_{n,\alpha} \|u_1\|_\alpha \dots \|u_n\|_\alpha
\en
with the same constant as before (so that $c_{n,\alpha} < 3^n$).

We use the latter bound to derive the second inequality (5.2). Let us split the density
of $X$ as $p = \varphi + \varphi^{1/\beta} v$, such that 
$\|v\|_\alpha^\alpha = \chi_\alpha(X,Z)$. Hence we get a decomposition
\bee
\hskip10mm \E\,e^{\frac{1}{2\beta} Z_n^2}
 & = &
\int_{\R^n} e^{\frac{n}{2\beta}\,\bar x^2}\,p(x_1)\dots p(x_n)\,dx_1 \dots dx_n \\
 & \hskip-45mm = & \hskip-24mm
\sum_{k=0}^n \frac{n!}{k!\,(n-k)!}\, \int_{\R^n} e^{\frac{n}{2\beta}\,\bar x^2}
\varphi(x_1)\dots \varphi(x_k)\,\varphi^{1/\beta}(x_{k+1}) v(x_{k+1}) \dots 
\varphi^{1/\beta}(x_n) v(x_n)
dx_1 \dots dx_n.
\ene
We apply (5.5) with $p_1$ to $p_k$ replaced by $\varphi$, and with $p_{k+1}$ to $p_n$
replaced with $\varphi^{1/\beta} v$
(that is, $u_j = \varphi^{1/\alpha}$ for $j \leq k$ and $u_j = v$ for $j > k$).
Moving the first term with $k = 0$ of this decomposition to the left, we then get the bound
$$
\big|\E\,e^{\frac{1}{2\beta} Z_n^2} - \E\,e^{\frac{1}{2\beta} Z^2}\big|  \, \leq \,
c_{n,\alpha} \sum_{k=1}^n \frac{n!}{k!\,(n-k)!}\, \|\varphi^{1/\alpha}\|_\alpha^k \, \|v\|_\alpha^{n-k} =
c_{n,\alpha} \big((1 + \|v\|_\alpha^n) - 1\big).
$$
\qed

\vskip10mm
\section{{\bf Connections with Fourier Transform}}
\setcounter{equation}{0}

\vskip2mm
\noindent
In the next sections, we restrict ourselves to the particular interesting index
$\alpha = 2$, that is, to the $\chi^2$-distance from the standard normal law,
$$
\chi^2(X,Z) = \int_{-\infty}^\infty \frac{p(x)^2}{\varphi(x)}\,dx - 1.
$$
In this case, necessary and sufficient conditions for finiteness of this divergence
may be given in terms of the characteristic function
$$
f(t) = \E\,e^{itX}, \qquad t \in \R.
$$

\vskip5mm
{\bf Proposition 6.1.} {\it The condition $\chi^2(X,Z) < \infty$ insures that $f(t)$ has 
square integrable derivatives of any order. Moreover, in that case
$$
1 + \chi^2(X,Z) \, = \, \frac{1}{\sqrt{2\pi}} \,
\sum_{n=0}^\infty \, \frac{1}{n!}\, \int_{-\infty}^\infty |f^{(n)}(t)|^2\,dt.
$$
}

\vskip2mm
{\bf Proof.} By the very definition,
$$
1 + \chi^2(X,Z) \, = \, \sqrt{2\pi}\, \sum_{n=0}^\infty \, \frac{1}{n!}\, 
\int_{-\infty}^\infty x^{2n} p(x)^2\,dx.
$$
We know that $f$ has finite derivatives of any order given by
$$
f^{(n)}(t) = \E\,(iX)^n\, e^{itX} = 
\int_{-\infty}^\infty (ix)^n e^{itx}\,p(x)\,dx.
$$
It remains to apply Plancherel's theorem.
\qed

\vskip5mm
In view of Proposition 4.1, existence of $\chi^2(X,Z)$ does not guarantee 
existence of the ``Gaussian" moment 
$\E\,e^{X^2/4}$. Nevertheless, it is true for the normalized convolution of the distribution 
of $X$ with itself, as indicated in Proposition 5.1.
In fact, in this case inequality (5.1) can be stated more precisely as
$$
\E\,e^{\frac{1}{4}\,(\frac{X + \widetilde X}{\sqrt{2}})^2} \leq 2\,(1+\chi^2(X,Z)),
$$
where $\widetilde X$ is an independent copy of $X$. Equivalently, there is a corresponding 
refinement of inequality (5.3) in Corollary 5.2 (without any convolution).

\vskip5mm
{\bf Proposition 6.2.} {\it In any case
$$
\frac {1}{\sqrt{2\pi}} \int_{-\infty}^{\infty} f(iy)^2\,e^{-2y^2}\,dy \le 1+\chi^2(X,Z).
$$
}

\vskip2mm
The argument is based on the following general observation which may be of independent interest.

\vskip5mm
{\bf Lemma 6.3.} {\it Given a function $p$ on the real line, suppose that
the function $g(x) = p(x)\,e^{x^2/4}$ belongs to $L^2$. Then the Fourier transforms
$$
f(t) = \int_{-\infty}^{\infty} e^{itx} p(x)dx, \qquad
\rho(t) = \int_{-\infty}^{\infty} e^{itx} g(x)\,dx
$$
are connected by the identity
\be
f(t) = \frac{1}{\sqrt{\pi}} \int_{-\infty}^{\infty} e^{-(t-u)^2}\rho(u)\,du \qquad (t \in \mathbb R),
\en
which may analytically be extended to the complex plane. Moreover, 
\be
\int_{-\infty}^{\infty} |f(iy)|^2\,e^{-2y^2}\,dy \, = \,
 \int_{-\infty}^{\infty} |\rho(t)|^2\,e^{-2t^2}\,dt.
\en
}

\vskip2mm
Thus, the characteristic function $f$ appears as the Weierstrass transform
of the function $g$. While Proposition 5.1 and its Corollary 5.2 are key ingredients of the proof
of Theorem 1.2, Lemma 6.2 can be used as an alternative approach to Theorem 1.1
for the particular case $\alpha = 2$. Lemma 6.3 and Proposition 6.2 can be adapted to cover 
the range $1 < \alpha \leq 2$ by considering the Fourier transform on the Lebesgue space $L^\alpha$.
However, these results do not extend to indexes $\alpha>2$.

Returning to the $L^2$-case, note that $g$ does not need to be integrable, so, one should
understand $\rho$ as a $L^2$-limit
$
\rho(t) = \lim_{N \rightarrow \infty} \int_{-N}^N e^{itx} g(x)\,dx
$
in the norm of the space $L^2$.

Note also that the second integral in (6.2) can be bounded by the squared
$L^2$-norm of $\rho$, which is, by the Plancherel theorem, equal to
$$
2\pi\, \|g\|_2^2
 \, = \, 2\pi\, \int_{-\infty}^{\infty} |p(x)|^2\,e^{x^2/2}\,dx.
$$
When $p$ is density of $X$, the last expression is exactly $\sqrt{2\pi}\,(1 + \chi^2(X,Z))$,
thus proving Proposition 6.2.

\vskip5mm
{\bf Proof of Lemma 6.3.} First assume that $p$ is a compactly supported;
in particular, both $p$ and $g$ are integrable and have analytic Fourier transforms.
By Fubini's theorem,
\bee
f(t)
 & = &
\int_{-\infty}^\infty e^{itx} g(x)\,
\bigg[\frac 1{\sqrt{\pi}} \int_{-\infty}^{\infty} e^{-ixu-u^2}\,du\bigg]\,dx \\
 &= &
\frac{1}{\sqrt{\pi}} \int_{-\infty}^\infty e^{-u^2}
\bigg[\int_{-\infty}^\infty e^{i(t-u)x} g(x)\,dx\bigg]\,du\\
 & = &
\frac{1}{\sqrt{\pi}} \int_{-\infty}^{\infty} e^{-(u-t)^2}\rho(u)\,du \ = \
\frac{1}{\sqrt{\pi}}\,e^{-t^2} \int_{-\infty}^{\infty}e^{2ut - u^2} \rho(u)\,du,
\ene
and we obtain (6.1). Moreover, a change of variable, we have
$$
\sqrt{\pi}\,f\Big(\frac{iz}{2}\Big)\,e^{-z^2/4} = 
\int_{-\infty}^{\infty}e^{izu - u^2} \rho(u)\,du \qquad (z \in \mathbb R),
$$
which means that the left-hand side represents the Fourier transform of the function
$e^{-u^2} \rho(u)$. Hence, by Plancherel's theorem,
\be
\|e^{-u^2} \rho(u)\|_2^2 = 
\frac{1}{2}\,\int_{-\infty}^\infty\,\Big|f\Big(\frac{iz}{2}\Big)\Big|^2\,e^{-z^2/4}\,dz =
\int_{-\infty}^{\infty} |f(iy)|^2\,e^{-2y^2} dy,
\en
thus proving (6.2).

In the general case, we have $p \in L^1 \cap L^2$, and arguing as 
in the proof of Proposition 4.1 (for the case $\alpha = 2$), we also get
$$
\int_{-\infty}^\infty
e^{cx^2}\,|p(x)|\,dx \leq \frac{C}{(1-4c)^{1/4}} < \infty \quad {\rm for \ all} \ \ 
c < \frac{1}{4},
$$
where $C^2 = \int_{-\infty}^\infty \frac{p(x)^2}{\varphi(x)}\,dx$.
In particular, $f$ is an entire function. Let $p_N$ be the restriction of $p$
to $[-N,N]$, $g_N(x) = p_N(x)\,e^{x^2/4}$, and put
$$
f_N(t) = \int_{-\infty}^{\infty} e^{itx} p_N(x)dx, \qquad
\rho_N(t) = \int_{-\infty}^{\infty} e^{itx} g_N(x)\,dx.
$$
According to the previous step, for all $t \in \mathbb R$,
\be
f_N(t) = \frac{1}{\sqrt{\pi}} \int_{-\infty}^{\infty} e^{-(t-u)^2}\rho_N(u)\,du.
\en
By the Lebesgue dominated convergence theorem, we have $f_N(t) \rightarrow f(t)$
for all real $t$ and $\|g_N - g\|_2 \rightarrow 0$ as $N \rightarrow \infty$.
By the continuity of the Fourier transform on $L^2$, we obtain
$\|\rho_N - \rho\|_2 \rightarrow 0$, which in turn implies
$$
\int_{-\infty}^{\infty} e^{-(t-u)^2}\rho_N(u)\,du \rightarrow 
\int_{-\infty}^{\infty} e^{-(t-u)^2}\rho(u)\,du.
$$
Hence, in the limit (6.4) yields the desired identity (6.1). Its right-hand side is well-defined
and finite for all complex $t$, and clearly represents an entire function.
Moreover, as before, one may apply Plancherel's theorem, leading to (6.3)
and therefore to (6.2).
\qed

\vskip5mm
\section{{\bf Exponential Series}}
\setcounter{equation}{0}

\vskip2mm
\noindent
The $\chi^2$-distance from the standard normal law on the real line admits a nice 
description in terms of a so-called exponential series (following Cram\'er's terminology)
as well. Let us some introduce basic notations and recall several well-known facts.
By $H_k$ we denote the $k$-th Chebyshev-Hermite polynomial
$$
H_k(x) = (-1)^k\, \big(e^{-x^2/2}\big)^{(k)}\, e^{x^2/2}, \qquad 
k = 0,1,2,\dots \ \ (x \in \R),
$$
so that $\varphi^{(k)}(x)\, = \,(-1)^k\, H_k(x)\varphi(x)$
in terms of the standard normal density. In particular,
\bee
H_0(x) = 1, & \quad H_2(x) = x^2 - 1, \\
H_1(x) = x, & \quad \ H_3(x) = x^3 - 3x.
\ene
Each $H_k$ is a polynomial of degree $k$ with integer coefficients, with leading coefficient 1. 
Depending on $k$ being even or odd, $H_k$ contains even resp. odd powers only.
These polynomials may be defined explicitly via
$$
H_k(x) = \E\,(x+iZ)^k, \quad Z \sim N(0,1).
$$
All $H_k$ are orthogonal to each other on the real line with weight function $\varphi(x)$, 
and moreover -- they form a complete orthogonal system in the Hilbert space 
$L^2(\R,\varphi(x)dx)$. Their $L^2$-norms are given by
$$
\E\, H_k(Z)^2 = \int_{-\infty}^\infty H_k(x)^2\, \varphi(x)\,dx = k!
$$
Equivalently, the Hermite functions $\varphi_k = H_k \varphi$
form a complete orthogonal system in $L^2(\R,\frac{dx}{\varphi(x)})$, 
and their $L^2$-norms in this space are given by
$
\int_{-\infty}^\infty \frac{\varphi_k(x)^2}{\varphi(x)}\,dx = k!
$
Summarizing we have:

\vskip5mm
{\bf Proposition 7.1.} {\it Any complex valued function 
$u = u(x)$ with
$
\int_{-\infty}^\infty |u(x)|^2\, e^{x^2/2}\,dx < \infty
$
admits a unique representation in the form of the orthogonal series
\be
u(x) = \varphi(x) \sum_{k=0}^\infty \frac{c_k}{k!}\, H_k(x),
\en
which converges in $L^2(\R,\frac{dx}{\varphi(x)})$.
The coefficients are given by
$
c_k = \int_{-\infty}^\infty u(x)\,H_k(x)\, dx,
$
and we have Parseval's identity
$$
\sum_{k=0}^\infty \frac{|c_k|^2}{k!} =
\int_{-\infty}^\infty \frac{|u(x)|^2}{\varphi(x)}\,dx.
$$
}

\vskip2mm
The functional series (7.1) representing $u$ is called an exponential series. 
The question of its pointwise convergence is rather delicate similar to 
the pointwise convergence of ordinary Fourier series based on trigonometric functions. 
In Cram\'er's paper [Cr], the following two propositions 
are stated, together with an explanation of the basic ingredients of the proof.

\vskip5mm
{\bf Proposition 7.2.} {\it If $u(x)$ is vanishing at infinity and
has a continuous derivative such that
$$
\int_{-\infty}^\infty |u'(x)|^2\, e^{x^2/2}\,dx < \infty,
$$
it may be developed in an exponential series, which is absolutely and 
uniformly convergent for $-\infty < x < \infty$.
}

\vskip5mm
{\bf Proposition 7.3.} {\it If $u(x)$ has bounded variation in every
finite interval, and if
$$
\int_{-\infty}^\infty |u(x)|\, e^{x^2/4}\,dx < \infty,
$$
then the exponential series for $u(x)$ converges to
$\frac{1}{2}\,(u(x+) + u(x-))$.
The convergence is uniform in every finite interval of continuity.
}

\vskip5mm
The integral condition of Proposition 7.3 is illustrated in [Cr] on 
the example of the Gaussian functions $u(x) = e^{-\lambda x^2}$ ($\lambda > 0$).
In this case, the corresponding exponential series can be explicitly 
computed, and at $x=0$ it is given by the series
$$
\frac{1}{\sqrt{2\lambda}} \ \sum_{k=0}^\infty \ 
\frac{(2k)!}{(k!)^2\,4^k}\, \Big(1 - \frac{1}{2\lambda}\Big)^k.
$$
This series is absolutely convergent for $\lambda>\frac{1}{4}$, 
simply convergent for $\lambda = \frac{1}{4}$ and divergent for
$\lambda < \frac{1}{4}$.

\vskip10mm
\section{{\bf Normal Moments}}
\setcounter{equation}{0}

\vskip2mm
\noindent
Let $X$ be a random variable with density $p$, and let $Z$ be a standard normal 
random variable (which is assumed to be independent of $X$). 
Applying Proposition 7.1 to $p$, we obtain the following: If
\be
\int_{-\infty}^\infty p(x)^2\, e^{x^2/2}\,dx < \infty,
\en
then $p$ admits a unique representation in the form of the exponential series
\be
p(x) = \varphi(x) \sum_{k=0}^\infty \frac{c_k}{k!}\, H_k(x),
\en
which converges in $L^2(\R,\frac{dx}{\varphi(x)})$. Here, 
the coefficients are given by
$$
c_k = \int_{-\infty}^\infty H_k(x)\,p(x)\, dx = \E H_k(X)= \E\,(X + iZ)^k,
$$
which we call the normal moments of $X$. In particular, $c_0 = 1$, $c_1 = \E X$,
$c_2 = \E X^2 - 1$. 

In general, these moments exist, as long as the $k$-th absolute moments of $X$ are 
finite. These moments are needed to develop the characteristic function of $X$
in a Taylor series around zero as follows:
\be
f(t) = \E\,e^{itX} = e^{-t^2/2}\, \sum_{k=0}^N \frac{c_k}{k!}\,(it)^k + o(|t|^N), \qquad
t \rightarrow 0.
\en
In particular, $c_k = 0$ for $k \geq 1$ in case $X$ is standard normal,
similarly to the property of the cumulants 
$$
\gamma_k(X) = \frac{d^k}{i^k\,dt^k}\, \log f(t)|_{t=0}
$$
with $k \geq 3$ (where we use the branch of the logarithm determined by $\log f(0) = 0$).

Let us emphasize one simple algebraic property of normal moments.

\vskip5mm
{\bf Proposition 8.1.} {\it Let $X$ be a random variable such that $\E X = 0$, 
$\E X^2 = 1$ and $\E\,|X|^k < \infty$ for some integer $k \geq 3$, and let
$Z \sim N(0,1)$. The following three properties are equivalent:

\vskip2mm
$a)$ \ $\gamma_j(X) = 0$\, for all\, $j = 3,\dots, k-1$;

\vskip1mm
$b)$ \ $\E H_j(X) = 0$ for all\, $j = 3,\dots, k-1$;

\vskip1mm
$c)$ \ $\E X^j = \E Z^j$ for all\, $j = 3,\dots, k-1$.

\vskip2mm
\noindent
In this case
\be
\gamma_k(X) \, = \, \E H_k(X) \, = \, \E X^k - \E Z^k.
\en
}

{\bf Proof.} Let us first describe the structure of the coefficients in (8.3) used for $N=k$. 
Repeated differentiation of the identity $f(t)\, e^{t^2/2} = \E\,e^{it(X + iZ)}$
yields $\frac{d^j}{i^j\,dt^j}\, \big[f(t)\, e^{t^2/2}\big]\big|_{t=0} = \E\,(X + iZ)^j$. 
Hence, we get indeed $c_j = \E H_j(X)$ for all $j \leq k$.

Now, assuming that $b)$ holds, the expansion (8.3) simplifes to 
\be
f(t) = e^{-t^2/2}\, \Big(1 + \frac{c_k}{k!}\,(it)^k\Big) + o(|t|^k), 
\en
so that
$\log f(t) = -\frac{1}{2}\,t^2 + \frac{c_k}{k!}\,(it)^k + o(|t|^k)$.
The latter expansion immediately yields $a)$. The argument may easily be reversed
in order to show that $a) \Rightarrow b)$ as well. Next, differentiating (8.5) $j$ times at zero, 
$j \leq k-1$, we get that $\E X^j = H_j(0)$. But, we obtain a similar equality
$\E Z^j = H_j(0)$ when writing (8.5) for $g(t) = e^{-t^2/2}$. Hence, $c)$ follows from $b)$.
Moreover, differentiating (8.5) $k$ times at zero, we arrive at $\E X^k = \E Z^k + c_k$,
which is the second equality in (8.4).
Again, the argument may be reversed in the sense that, starting from $c)$, we
obtain (8.5) and therefore $b)$. Thus, all the three properties are equivalent.

Finally, the first equality in (8.4) is obtained when differentiating the expression
$\log f(t) = -\frac{1}{2}\,t^2 + \frac{c_k}{k!}\,(it)^k + o(|t|^k)$ $k$ times.
\qed

\vskip5mm
In general (without the above conditions $a)-c)$), the moments of $X$ may be expressed 
easily in terms of the normal moments. Indeed, the Chebyshev-Hermite polynomials 
have generating function
$$
\sum_{k=0}^\infty H_k(x)\, \frac{z^k}{k!} = e^{xz - z^2/2},
$$
which follows, for example, from the identity $H_k(x) = \E\,(x+iZ)^k$. Here 
$x$ may $z$ may be any complex numbers. Equivalently,
$$
e^{xz} \, = \,
e^{z^2/2} \sum_{i=0}^\infty H_i(x)\, \frac{z^i}{i!} \, = \,
\sum_{i,j=0}^\infty H_i(x)\,\frac{z^{i+2j}}{i!j!\, 2^j}.
$$
Expanding $e^{xz}$ into the power series and comparing
the coefficients in front of $z^k$, we get
$$
x^k = k! \sum_{j=0}^{[k/2]} \frac{1}{(k-2j)!\, j!\, 2^j}\, H_{k-2j}(x).
$$
Hence, if $\E\,|X|^k < \infty$, then
\be
\E X^k = k! \sum_{j=0}^{[k/2]} \frac{1}{(k-2j)!\, j!\, 2^j}\, 
\E H_{k-2j}(X).
\en

Now, let us describe the connection between normal moments and the $\chi^2$-distance.
The series in (8.3) is absolutely convergent as $N \rightarrow \infty$, when $f$ is analytic
in the complex plane. Hence we have the expansion
\be
f(t) = e^{-t^2/2}\, \sum_{k=0}^\infty \frac{c_k}{k!}\,(it)^k, \qquad t \in \C,
\en
which holds, in particular, assuming condition (8.1). Moreover, using the Parseval identity as in
Proposition 7.1, we have 
\be
\sum_{k=0}^\infty \frac{|c_k|^2}{k!} =
\int_{-\infty}^\infty \frac{p(x)^2}{\varphi(x)}\,dx.
\en
Since right-hand side is related to $\chi^2$-distance from the standard normal law, 
we arrive at the following relation:

\vskip5mm
{\bf Proposition 8.2.} {\it If $\chi^2(X,Z) < \infty$, then
\be
\chi^2(X,Z) \, = \, \sum_{k=1}^\infty \frac{1}{k!}\,(\E H_k(X))^2.
\en
}

For the quadratic Renyi divergence, we thus have
$$
D_2(X||Z) = \log(1 + \chi^2(X,Z)) = \log \,
\sum_{k=0}^\infty \frac{1}{k!}\,(\E H_k(X))^2.
$$

Recall that, if $\chi^2(X,Z) < \infty$, then $X$ has finite moments of any order, 
and moreover, $\E\,e^{cX^2} < \infty$ for any $c < \frac{1}{4}$. Hence, the 
normal moments $\E H_k(X)$ are well defined and finite, so that the representation 
for $\chi^2(X,Z)$ makes sense.

We now show a converse to Proposition 8.2.

\vskip5mm
{\bf Proposition 8.3.} {\it Let $X$ be a random variable with finite
moments of any order. If the series in $(8.9)$ is convergent,
then $X$ has an absolutely continuous distribution with finite distance $\chi^2(X,Z)$.
}

\vskip5mm
It looks surprising that a simple sufficient condition for the existence of a density 
$p$ of $X$ can be formulated in terms of moments of $X$, only. 
Note that if $X$ is bounded, then it has finite moments of any order, and 
the property $\chi^2(X,Z) < \infty$ just means that $p$ is in $L^2$.
Thus, we have:

\vskip5mm
{\bf Corollary 8.4.} {\it A bounded random variable $X$ has an absolutely
continuous distribution with a square integrable density, if and only if
the series in $(8.9)$ is convergent.
}

\vskip5mm
{\bf Proof of Proposition 8.3.}
Let $C^2 = \sum_{k=0}^\infty \, \frac{1}{k!}\,(\E H_k(X))^2$
be finite ($C \geq 1$). Then $|\E H_k(X)| \leq C\sqrt{k!}$ and from the formula 
(8.6) we get
\bee
|\E X^k| 
 & \leq &
k! \sum_{j=0}^{[k/2]} \frac{1}{(k-2j)!\,j!\, 2^j}\, |\E H_{k-2j}(X)|
 \ \leq \ 
C k! \sum_{j=0}^{[k/2]} \frac{1}{\sqrt{(k-2j)!}\, j!\, 2^j}.
\ene
In particular,
$$
\E X^{2k} \, \leq \, C (2k)!\, \sum_{j=0}^{k} \frac{1}{\sqrt{(2k-2j)!}\ j!\, 2^j}.
$$
Using $\frac{(2k)!}{(2k-2j)!} \leq (2k)^{2j}$, we obtain that
\bee
\E X^{2k} 
 & \leq & 
C \sqrt{(2k)!} \ 
\sum_{j=0}^{k} \frac{\sqrt{(2k)!}}{\sqrt{(2k-2j)!}\ j!\, 2^j} \\
 & \leq & 
C \sqrt{(2k)!}\ \sum_{j=0}^{k} \frac{(2k)^j}{j!\, 2^j} \ < \ 
C \sqrt{(2k)!}\ \sum_{j=0}^\infty \frac{k^j}{j!} \ = \
Ce^k \sqrt{(2k)!}
\ene

Thus, $\E X^{2k} < Ce^k \sqrt{(2k)!}$ for all $k$.
This estimate implies that
$\E\,e^{cX^2} < \infty$ for some $c>0$. In particular, $X$ has an entire 
characteristic function $f(t) = \E\,e^{itX}$ which thus admits a power series 
representation (8.7), where necessarily $c_k = \E H_k(X)$.
Consider the $N$-th partial sum of that series,
$$
f_N(t) = e^{-t^2/2} \sum_{k=0}^N c_k\,\frac{(it)^k}{k!}.
$$
It represents the Fourier transform of the function
$p_N(x) = \varphi(x) \sum_{k=0}^N c_k\,\frac{H_k(x)}{k!}$
which is the $N$-th partial sum of the exponential series in (8.2). Since, by the assumption,
$$
\sum_{k=0}^\infty \frac{c_k^2}{k!} < \infty,
$$
$p_N$ converge to some $p$ in $L^2(\R,\frac{dx}{\varphi(x)})$, by Proposition 7.1. 
In particular, $p_N$ converge in 
$L^2(\R,dx)$, and by Plancherel's theorem, $f_N$ also converge 
in $L^2(\R,dx)$ to the Fourier transform $\hat p$ of $p$. But $f_N(t) \rightarrow f(t)$
for all $t$, so $f(t) = \hat p(t)$ almost everywhere.
Thus we conclude that $f$ belongs to $L^2(\R,dx)$ and is equal to
the Fourier transform of $p$. Hence, $X$ has an absolutely continuous
distribution, and $p$ is density of $X$.

It remains to use once more the orthogonal series (8.2). 
By Proposition 7.1, we have Parseval's equality (8.8), which means that
$
\chi^2(X,Z) = \sum_{k=0}^\infty \frac{c_k^2}{k!} < \infty.
$
\qed

\vskip5mm
There is a natural generalization of the identity (8.9) in terms of
the $\chi^2$-distance for the random variables
$$
X_t = \sqrt{t}\, X + \sqrt{1-t}\, Z, \qquad 0 \leq t \leq 1,
$$
where $Z \sim N(0,1)$ is independent of $X$.

\vskip5mm
{\bf Proposition 8.5.} {\it If $\chi^2(X,Z) < \infty$, then, for all $t \in [0,1]$,
$$
\chi^2(X_t,Z) \, = \, \sum_{k=1}^\infty \frac{t^k}{k!}\,(\E H_k(X))^2.
$$
}

This in turn yields another description of the normal moments via the derivatives
of the $\chi^2$-distance:
$$
(\E H_k(X))^2 \, = \, \frac{d^k t}{dt^k}\,\chi^2(X_t,Z)\big|_{t=0}, \quad k = 1,2,\dots
$$

\vskip2mm
{\bf Proof.}
It is known, e.g., as a direct consequence of the identity $H_k(x) = \E\,(x+iZ)^k$, 
that the Hermite polynomials satisfy the binomial formula
\be
H_k(ax+by) = \sum_{i=0}^k C_k^i\, a^i b^{k-i}\, H_i(x) H_{k-i}(y), \qquad
x,y \in \R, 
\en
whenever $a^2 + b^2 = 1$. In particular,
$
\E H_k(aX+bZ) = a^k\,\E H_k(X),
$
which may be used in the formula (8.9) with $a = \sqrt{t}$ and $b = \sqrt{1-t}$.
\qed

\vskip10mm
\section{{\bf Behavior of R\'enyi divergence under Convolutions}}
\setcounter{equation}{0}

\vskip2mm
\noindent
The obvious question, when describing convergence in the central limit theorem
in the $D_\alpha$-distance is, it remain finite for sums of independent summands with 
finite $D_\alpha$-distances? The answer is affirmative and is made precise in the following:

\vskip5mm
{\bf Proposition 9.1.} {\it Let $X$ and $Y$ be independent random variables.
Given $\alpha > 1$, for all $a,b \in \R$ such that $a^2 + b^2 = 1$, we have
$$
D_\alpha(aX+bY || Z) \leq D_\alpha(X||Z) + D_\alpha(Y||Z),
$$
where $Z \sim N(0,1)$. Equivalently,
\be
1 + (\alpha - 1)\,T_\alpha(aX+bY||Z) \, \leq \, 
\big(1+(\alpha - 1)\,T_\alpha(X||Z)\big)\, \big(1 + (\alpha - 1)\,T_\alpha(Y||Z)\big).
\en
}

\vskip2mm
The statement may be extended by induction to finitely many 
independent summands $X_1,\dots,X_n$ by the relation
$$
D_\alpha(a_1X_1+ \dots + a_n X_n || Z) \leq D_\alpha(X_1||Z) + \dots + D_\alpha(X_n||Z),
$$
where $a_1^2 + \dots + a_n^2 = 1$.
Note that for the relative entropy ($\alpha=1$), there is a much stronger property, namely
$$
D(a_1X_1+ \dots + a_n X_n || Z) \leq \max\{D(X_1||Z),\dots,D(X_n||Z)\},
$$
which follows from the entropy power inequality (cf. [D-C-T]).
However, this is no longer true for $D_\alpha$.
Nevertheless, for the normalized sums $Z_n = (X_1 + \dots + X_n)/\sqrt{n}$
with i.i.d. summands, Proposition 9.1 guarantees a sublinear growth of the 
R\'enyi divergence with respect to $n$. More precisely, we have
\begin{equation}
D_\alpha(Z_n||Z) \leq n D_\alpha(X_1||Z).
\end{equation}

\vskip2mm
{\bf Proof of Proposition 9.1.} Let $Z$ be an independent copy of $Z$, so that
the random vector $\widetilde Z = (Z,Z')$ is standard normal in $\R^2$.
From Definition 2.1 it follows that the R\'enyi distance of the random vector 
$\widetilde X = (X,Y)$ to $\widetilde Z$ is given by
$$
D_\alpha(\widetilde X||\widetilde Z) = D_\alpha(X||Z) + D_\alpha(Y||Z').
$$
Hence, by the contractivity property (2.1), cf. Proposition 2.3, we get
$$
D_\alpha\big(S(\widetilde X)||S(\widetilde Z)\big) \leq
D_\alpha(X||Z) + D_\alpha(Y||Z'),
$$
for any Borel measurable function $S:\R^2 \rightarrow \R$. It remains to apply
this inequality to the linear function $S(x,y) =ax + by$.
\qed

\vskip5mm
Let us describe a simple alternative argument in the case $\alpha = 2$, which relies 
upon normal moments only. One may assume that both $D_2(X||Z)$ and $D_2(Y||Z)$ 
are finite, so that  $X$ and $Y$ have finite moments of any order. 
In addition, without loss of generality, let $a,b>0$.

From the binomial formula (8.10) it follows that
$$
\E H_k(aX+bY) = \sum_{i=0}^k C_k^i\, a^i b^{k-i}\, 
\E H_i(X)\, \E H_{k-i}(Y).
$$
By Cauchy's inequality,
\bee
(\E H_k(aX+bY))^2 
 & \leq &
\sum_{i=0}^k C_k^i\, (a^i b^{k-i})^2\, 
\sum_{i=0}^k C_k^i\, (\E H_i(X))^2\, (\E H_{k-i}(Y))^2 \\
 & = &
\sum_{i=0}^k C_k^i\, (\E H_i(X))^2\, (\E H_{k-i}(Y))^2.
\ene
This gives
$$
\frac{(\E H_k(aX+bY))^2}{k!} \, \leq \,
\sum_{i=0}^k \frac{(\E H_i(X))^2}{i!}\, \frac{(\E H_{k-i}(Y))^2}{(k-i)!},
$$
and summation over all integers $k \geq 0$ leads to
$$
\sum_{k=0}^\infty \frac{(\E H_k(aX+bY))^2}{k!} \, \leq \,
\sum_{i=0}^\infty \frac{(\E H_i(X))^2}{i!}\, 
\sum_{j=0}^\infty\frac{(\E H_j(Y))^2}{j!}.
$$
But, by Proposition 8.2, this inequality is the same as 
$$
1 + \chi^2(aX+bY,Z) \, \leq \, \big(1+\chi^2(X,Z)\big)\, \big(1 + \chi^2(Y,Z)\big),
$$
which is exactly (9.1) for $\alpha=2$.

One may also ask whether or not $\chi^2(aX + bY,Z)$ remains finite, when
$\chi^2(X,Z)$ is finite, and $Y$ is ``small" enough. If $p$ is density of $X$, the
density of $aX+bY$ is given by
$$
q(x) = \frac{1}{|a|}\,\E\,p\Big(\frac{x - bY}{a}\Big), \qquad x \in \R,
$$
which is a convex mixture of densities on the line. Applying Cauchy's inequality,
we have
$$
\frac{q(x)^2}{\varphi(x)} \, \leq \, \frac{1}{a^2}\,
\E\,\frac{p(\frac{x - bY}{a})^2}{\varphi(x)},
$$
and using $(ax + by)^2 \leq x^2 + y^2$, we get an elementary bound
\bee
\int_{-\infty}^\infty \frac{q(x)^2}{\varphi(x)}\,dx
  & \leq & 
\frac{1}{|a|}\, \E\, \int_{-\infty}^\infty \frac{p(x)^2}{\varphi(ax + bY)}\,dx \\
 & \leq &
\frac{1}{|a|}\, \E\, \int_{-\infty}^\infty 
\sqrt{2\pi}\, p(x)^2\, e^{\frac{1}{2} (x^2 + Y^2)}\,dx \ = \ 
\frac{1}{|a|}\,(1 + \chi^2(X,Z))\,\E\,e^{Y^2/2}.
\ene
That is, we arrive at:

\vskip5mm
{\bf Proposition 9.2.} {\it Let $X$ and $Y$ be independent random variables.
For all $a,b \in \R$ such that $a^2 + b^2 = 1$, we have
$$
1 + \chi^2(aX+bY,Z) \, \leq \, \frac{1}{|a|}\,(1 + \chi^2(X,Z))\,\E\,e^{Y^2/2},
\qquad Z \sim N(0,1).
$$
}

\vskip2mm
Let us now describe two examples of i.i.d. random variables $X,X_1,\dots,X_n$
such that for the normalized sums $Z_n = (X_1 + \dots + X_n)/\sqrt{n}$, we have
\be
\chi^2(Z_1,Z) = \dots = \chi^2(Z_{n_0 - 1},Z) = \infty, \qquad
\chi^2(Z_{n_0},Z) < \infty,
\en
where $n_0 > 1$ is a given prescribed integer.

\vskip5mm
{\bf Example 9.3.} Suppose that $X$ has density of the form
\be
p(x) = \int_0^\infty \frac{1}{\sigma \sqrt{2\pi}}\,e^{-x^2/2\sigma^2}
\,d\pi(\sigma^2),  \qquad x \in \R,                                     
\en
where $\pi$ is a probability measure on the positive half-axis. The existence of
$\chi^2(X,Z)$ implies that $\sigma^2 < 2$ for $\pi$-almost 
all $\sigma^2$, i.e., $\pi$ should be supported on the interval $(0,2)$. Squaring
the equality (9.4) and integrating over $x$, we find that
$$
1 + \chi^2(X,Z) \, = \, \int_{-\infty}^\infty \frac{p(x)^2}{\varphi(x)}\,dx \, = \,
\int_0^2 \!\!\int_0^2 \frac{1}{\sqrt{\sigma_1^2 + \sigma_2^2 - \sigma_1^2\sigma_2^2}} 
\ d\pi(\sigma_1^2)\, d\pi(\sigma_2^2).
$$
It is easy to see that the last double integral is convergent, if and only if
$$
\int_0^1 \!\!\int_0^1 \frac{1}{\sqrt{\sigma_1^2 + \sigma_2^2}}
\ d\pi(\sigma_1^2)\, d\pi(\sigma_2^2) < \infty \ \ {\rm and} \ \ 
\int_1^2 \!\!\int_1^2 \frac{1}{\sqrt{4 - (\sigma_1^2 + \sigma_2^2)}} 
\ d\pi(\sigma_1^2)\, d\pi(\sigma_2^2) < \infty.
$$
These conditions may be simplified in terms of the distribution function
$F(\ep) = \pi\{\sigma^2 \leq \ep\}$, $0 \leq \ep \leq 2$, by noting that
$$
F(\ep/2)^2 \leq (\pi \otimes \pi)\{\sigma_1^2 + \sigma_2^2 \leq \ep\} \leq F(\ep)^2.
$$
Hence, the first integral is convergent, if and only if
$$
\int_0^1 \frac{1}{\sqrt{\ep}}\,dF(\ep)^2 \, = \, 
F(1-)^2 + \frac{1}{2} \int_0^1 \frac{F(\ep)^2}{\ep^{3/2}}\,d\ep
$$
is finite. A similar description applies to the second double integral.

Let us summarize: We have $\chi^2(X,Z) < \infty$ for the random variable $X$
with density (9.4), if and only if the mixing probability measure $\pi$ is supported 
on the interval $(0,2)$, and its distribution function $F$ satisfies
\be
\int_0^1 \frac{F(\ep)^2}{\ep^{3/2}}\,d\ep < \infty, \qquad
\int_1^2 \frac{(1 - F(\ep))^2}{(2 - \ep)^{3/2}}\,d\ep < \infty.
\en

\vskip2mm
Based on this description, we now investigate convolutions. Note that $Z_n$ has density 
of a similar type as before
$$
p_n(x) = \int_0^\infty \frac{1}{\sigma \sqrt{2\pi}}\,e^{-x^2/2\sigma^2}
\,d\pi_n(\sigma^2).                                 
$$
More precisely, if $\xi_1,\dots,\xi_n$ are independent copies of a random variable $\xi$ 
distributed according to $\pi$, then the mixing measure $\pi_n$ can be recognized as
the distribution of the normalized sum $S_n = \frac{1}{n}\,(\xi_1 + \dots + \xi_n)$.
Therefore, by (9.5), $\chi^2(Z_n,Z) < \infty$, if and only if $\P\{S_n < 2\} = 1$ 
(which is equivalent to the property that $\pi$ is supported on $(2,\infty)$), and
$$
\int_0^1 \frac{F_n(\ep)^2}{\ep^{3/2}}\,d\ep < \infty, \qquad
\int_1^2 \frac{(1 - F_n(\ep))^2}{(2 - \ep)^{3/2}}\,d\ep < \infty,
$$
where $F_n$ is the distribution function of $S_n$. Since
$
F(\ep/n)^n \leq F_n(\ep) \leq F(\ep)^n,
$
which is needed near zero, and using similar relations near the point 2,
these conditions may be simplified to
\be
\int_0^1 \frac{F(\ep)^{2n}}{\ep^{3/2}}\,d\ep < \infty, \qquad
\int_1^2 \frac{(1 - F(\ep))^{2n}}{(2 - \ep)^{3/2}}\,d\ep < \infty.
\en

Now, for simplicity, suppose that $\pi$ is supported on $(0,1)$, so that the second
integral in (9.6) is vanishing, and let $F(\ep) \sim \ep^\kappa$ for $\ep \rightarrow 0$ 
with parameter $\kappa > 0$ (where the equivalence is understood
up to a positive factor). Then, the first integral in (9.6) will be finite, if and only if 
$n > 1/(4\kappa)$. Choosing $\kappa = \frac{1}{4 (n_0 - 1)}$,
we obtain the required property (9.3).

\vskip5mm
{\bf Example 9.4.} Consider a density of the form
$$
p(x) = \frac{a_k}{1 + |x|^{1/2k}}\,e^{-x^2/4},  \qquad x \in \R,                                     
$$
where $a_k$ is a normalizing constant, $k = n_0 - 1$, and let $f_1$ denote its
Fourier transform (i.e., the characteristic function). Define the distribution of $X$ via its 
characteristic function
$$
f(t) = \alpha f_1(t) + (1 - \alpha)\,\frac{\sin(\gamma t)}{\gamma t}
$$
with a sufficiently small parameter $\alpha > 0$ and
$\gamma = \sqrt{3\,\frac{1 + \alpha f_1''(0)}{1 - \alpha}}$. It is easy to check that
$f''(0) = -1$, which guarantees that $\E X= 0$, $\E X^2 = 1$. Furthermore, it is 
not difficult to show that the densities $p_n$ of $Z_n$ admit the two-sided bounds
$$
\frac{b_n'}{1 + |x|^{n/2k}}\,e^{-x^2/4} \leq p_n(x) \leq
\frac{b_n''}{1 + |x|^{n/2k}}\,e^{-x^2/4} \qquad (x \in \R),
$$
up to some positive $n$-dependent factors. Hence, again we arrive at the property (9.3).

\vskip10mm
\section{{\bf Superadditivity of $\chi^2$ with Respect to Marginals}}
\setcounter{equation}{0}

\vskip2mm
\noindent
A multidimensional version of Theorem 1.1 requires to involve some other properties of 
the $\chi^2$-distance in higher dimensions. The contractivity under mappings, 
$$
\chi^2(S(X),S(Z)) \leq \chi^2(X,Z),
$$
has already been shown in Proposition 2.3 in a general setting.
This inequality may be considerably sharpened, when distance is measured to the 
standard normal law in $\Omega = \R^d$. In order to compare the behavior of
$\chi^2$-divergence with often used information-theoretic quantities, recall the definition 
of Shannon entropy and Fisher information,
$$
h(X) = -\int_{\R^d} p(x)\,\log p(x)\,dx, \qquad
I(X) = \int_{\R^d} \frac{|\nabla p(x)|^2}{p(x)}\,dx,
$$
where $X$ is a random vector in $\R^d$ with density $p$ (assuming that the above 
integrals are well-defined). These functionals are known to be subadditive and 
superadditive with respect to the components: If we write $X = (X',X'')$ with 
$X' \in \R^{d_1}$, $X'' \in \R^{d_2}$ ($d_1 + d_2 = d$), then one always has
\be
h(X) \leq  h(X') + h(X''), \qquad 
I(X) \geq I(X') + I(X'')
\en
cf. [L], [C]. Both $h(X)$ and $I(X)$ themselves are not yet distances,
so one also considers the relative entropy and the relative Fisher information with respect
to other distributions. In particular, in case of the standard normal random vector 
$Z \sim N(0,{\rm I}_d)$ and random vectors $X$ with mean zero and identity covariance 
matrix ${\rm I}_d$, they are given by
$$
D(X||Z) = h(Z) - h(X), \qquad I(X||Z) = I(X) - I(Z). 
$$
Hence, it immediately follows from (10.1) that these information-theoretic distances are 
both superadditive, that is,
$$
D(X||Z) \geq  D(X'||Z') + D(X''||Z''), \qquad 
I(X||Z) \geq I(X'||Z') + I(X''||Z''),
$$
where $Z'$ and $Z''$ are standard normal in $\R^{d_1}$ and $\R^{d_2}$ respectively
(both inequalities become equalities, when $X'$ and $X''$ are independent).

We now establish a similar property for the $\chi^2$-distance, which can be more
convenient stated in the setting of a Euclidean space $H$, say of dimension $d$, 
with norm $|\cdot|$ and inner product $\left<\cdot,\cdot\right>$. If $X$ is a random 
vector in $H$ with density $p$, and $Z \sim N(0,{\rm I}_d)$
is a normal random vector with mean zero and an identity covariance operator
${\rm I}_d$, then (according to the abstract definition),
$$
\chi^2(X,Z) = \int_H \, \frac{p(x)^2}{\varphi(x)}\, dx - 1 = 
\int_H \frac{(p(x) - \varphi_d(x))^2}{\varphi(x)}\, dx,
$$
where $\varphi(x) = \frac{1}{(2\pi)^{d/2}}\,e^{-|x|^2/2}$ ($x \in H$)
is the density of $Z$.

\vskip5mm
{\bf Proposition 10.1.} {\it Given a random vector $X$ in $H$ and an orthogonal decomposition 
$H = H' \oplus H''$ into two linear subspaces $H',H'' \subset H$ 
of dimensions $d_1,d_2 \geq 1$, for orthogonal projections 
$X' = {\rm Proj}_{H'}(X)$, $X'' = {\rm Proj}_{H''}(X)$, we have
\be
\chi^2(X,Z) \, \geq \, \chi^2(X',Z') + \chi^2(X'',Z''),
\en
where $Z,Z',Z''$ are standard normal random vectors in $H,H',H''$, respectively.
}

\vskip5mm
Note, however, that (10.2) won't become an equality for independent components
$X',X''$.

\vskip5mm
{\bf Proof.} 
Let $H = \R^d$ and $X = (\xi_1,\dots,\xi_d)$.
Note that $\chi^2(X,Z)$ is invariant under orthogonal transformations $U$ of the space, i.e.,
$\chi^2(U(X),Z) = \chi^2(X,Z)$. Hence, without loss of generality, one may assume that
$X' = (\xi_1,\dots,\xi_{d_1})$ and $X'' = (\xi_{d_1 +1},\dots,\xi_d)$. Moreover, to simplify 
the argument (notationally), let $d_1 = d_2 = 1$.

The finiteness of the distance $\chi^2(X,Z)$ means that the random vector $X = (\xi_1,\xi_2)$ 
has density $p = p(x_1,x_2)$ ($x_i \in \R$) such that
$$
\int_{-\infty}^\infty \!\int_{-\infty}^\infty p(x_1,x_2)^2\, e^{(x_1^2 + x_2^2)/2}\,dx_1dx_2 < \infty.
$$
Since the Hermite functions 
$\varphi_{k_1,k_2}(x_1,x_2) = \varphi(x_1)\varphi(x_2)\, H_{k_1}(x_1) H_{k_2}(x_2)$
form a complete orthogonal system in $L^2(\R^2)$ (where now $\varphi$ denotes the
one dimensional standard normal density), the density $p$
admits a unique representation in the form of the exponential series
\be
p(x_1,x_2) \, = \, \varphi(x_1)\varphi(x_2) \sum_{k_1=0}^\infty \sum_{k_2=0}^\infty \,
\frac{c_{k_1,k_2}}{k_1! k_2!}\, H_{k_1}(x_1) H_{k_2}(x_2),
\en
which converges in $L^2(\R,\frac{dx_1 dx_2}{\varphi(x_1)\varphi(x_2)})$. Here 
the coefficients (mutual normal moments) are given by
$$
c_{k_1,k_2} \, = \, \int_{-\infty}^\infty \!\int_{-\infty}^\infty 
H_{k_1}(x_1) H_{k_2}(x_2)\,p(x_1,x_2)\,dx_1dx_2 \, = \,
\E H_{k_1}(\xi_1) H_{k_2}(\xi_2),
$$
and we have Parseval's equality
\be
1 + \chi^2(X,Z) \, = \, \int_{-\infty}^\infty \!\int_{-\infty}^\infty
\frac{p(x_1,x_2)^2}{\varphi(x_1)\varphi(x_2)}\,dx_1dx_2 \ = 
\sum_{k_1=0}^\infty \sum_{k_2=0}^\infty\,\frac{c_{k_1,k_2}^2}{k_1! k_2!}.
\en

Now, integrating (10.3) over $x_2$ and separately over $x_1$, we obtain similar representations
for the marginal densities
$$
p_1(x_1) \, = \, \varphi(x_1) \sum_{k_1=0}^\infty\,
\frac{c_{k_1,0}}{k_1!}\, H_{k_1}(x_1), \qquad
p_2(x_2) \, = \, \varphi(x_2) \sum_{k_2=0}^\infty\,
\frac{c_{0,k_2}}{k_2!}\, H_{k_2}(x_2),
$$
hence, by Proposition 8.1,
$$
\chi^2(\xi_1,\xi) \, = \, \sum_{k_1 = 1}^\infty \frac{c_{k_1,0}^2}{k_1!}, \qquad
\chi^2(\xi_2,\xi) \, = \, \sum_{k_2 = 1}^\infty \frac{c_{0,k_2}^2}{k_2!} \qquad
(\xi \sim N(0,1)).
$$
Obviously, the quantities $\chi^2(\xi_1,\xi)$ and $\chi^2(\xi_2,\xi)$ 
appear as summands in (10.4), thus showing the inequality.
\qed

\newpage
\noindent
{\large Part II: The R\'enyi Divergence in the Central Limit Theorem}

\vskip15mm
\section{{\bf Asymptotic Expansions and Lower Bounds}}
\setcounter{equation}{0}

\noindent
Let $X,X_1,X_2,\dots$ be independent identically distributed
random variables such that $\E X = 0$, $\E X^2 = 1$, with characteristic
function $f(t) = \E\,e^{itX}$. Then the normalized sums
$$
Z_n = \frac{X_1 + \dots + X_n}{\sqrt n}
$$
weakly converge in distribution to the standard normal law:
$Z_n \Rightarrow Z$ for $Z \sim N(0,1)$. In this connection
the following question arises:
When is it true that $D_\alpha(Z_n||Z) \rightarrow~0$ or equivalently
$T_\alpha(Z_n||Z) \rightarrow 0$ as $n \rightarrow \infty$?
And if so, what is the rate of convergence?

We shall give a complete solution of this problem in the next sections. First we shall describe here 
asymptotic expansions for ``truncated" $T_\alpha$-distances, which yield
reasonable lower bounds for $T_\alpha(Z_n||Z)$. More precisely, given $M>0$,
we have an obvious estimate
\be
T_\alpha(Z_n||Z) \, \ge \, \frac{1}{\alpha - 1}\,(I(M) - 1)
\en
with
\be
I(M) = \int_{|x|\le M}\bigg(\frac{p_n(x)}{\varphi(x)}\bigg)^\alpha\,\varphi(x)\,dx,
\en
where $p_n$ denotes the density of $Z_n$.
We will see that, under suitable conditions (like the ones in Theorems 1.1-1.2), 
while choosing
$$
M = M_n(s) = \sqrt{2(s-1)\log n}
$$ 
with a fixed integer $s \geq 2$, inequality (11.1) can be reversed up to an error 
term of order $o(n^{-(s-1)})$.
This reduces our task to the study of the asymptotic behavior of the integrals 
$I(M_n(s))$, using the following result due to Petrov (cf. [Pe1-2], [B-C-G1]).

\vskip5mm
{\bf Proposition 11.1.}
{\it Suppose that $X$ has a finite absolute moment of order $k\ge 3$, and assume 
that $Z_n$ admits a density in $L^2$ for some $n$. Then, for all 
$n$ large enough, $Z_n$ have continuous bounded densities $p_n$ satisfying
uniformly in $-\infty < x < \infty$
\be
p_n(x)=\varphi(x) +  \varphi(x) \sum_{\nu=1}^{k-2}
\frac{q_{\nu}(x)}{n^{\nu/2}} + o\Big(\frac 1{n^{(k-2)/2}}\Big)\frac 1{1+|x|^k}.
\en
}

\vskip2mm
In this formula
\be
q_{\nu}(x) = \sum H_{\nu+2l}(x)
\prod_{m=1}^{\nu}\frac 1{k_m!}\Big(\frac{\gamma_{m+2}}{(m+2)!}\Big)^{k_m}, 
\en
where $\gamma_r$ denotes the $r$-th cumulant of $X$. 
The summation extends over all non-negative integer solutions $(k_1,k_2,\dots,k_{\nu})$ 
to the equation $k_1+2k_2+\dots+\nu k_{\nu}=\nu$, and where we put 
$l=k_1+k_2+\dots+k_{\nu}$. The sum in (11.3) defines a polynomial in $x$ 
of degree at most $3(k-2)$.

For example, for $k=3$ (11.3) yields
$$
p_n(x) = \varphi(x) + \frac{\gamma_3}{3!\sqrt{n}}\,H_3(x) \varphi(x)
 + o\Big(\frac{1}{\sqrt{n}}\Big)\frac 1{1+|x|^3},
$$
where $\gamma_3 = \alpha_3 = \E X^3$ and $H_3(x) = x^3 - 3x$.
More generally, if the first cumulants of $X$ up to order $k-1$ are vanishing 
(i.e., the fist $k-1$ moments of $X$ are 
the same as for a standard normal law), then (11.3) simplifies to the expression
$$
p_n(x)=\varphi(x) + \frac{\gamma_k}{k!}\,H_k(x) \varphi(x) \, n^{-\frac{k-2}{2}}
 + o\big(n^{-\frac{k-2}{2}}\big)\frac 1{1+|x|^k} \qquad
(\gamma_3 = \dots = \gamma_{k-1} = 0).
$$

This local limit theorem may be used to derive:

\vskip5mm
{\bf Lemma 11.2.}
{\it Under the assumptions of Proposition $11.1$ with $k = 2s$ $(s \geq 2)$, 
the following expansion holds
\be
I(M_n(s)) = 1 + \sum_{j=1}^{s-1} \frac{b_j}{n^j} + o\big(n^{-(s-1)}\big)
\en
with
\be
b_j = \sum \frac{\alpha (\alpha - 1) \dots (\alpha - m + 1)}{m_1! \dots m_{2j-1}!}\,
\int_{-\infty}^\infty q_1(x)^{m_1} \dots q_{2j-1}(x)^{m_{2j-1}}\,\varphi(x)\,dx.
\en
Here the sum extends over all non-negative integer numbers $m_1, \dots, m_{2j-1}$ 
such that $m_1 + 2m_2 + \dots + (2j-1)\,m_{2j-1} = 2j$, and where 
$m = m_1 + \dots + m_{2j-1}$. In particular,
if $\gamma_j=0$ for $j=3,\dots,s-1$, $s \geq 3$, then
\be
I(M_n(s)) = 1 + \alpha (\alpha-1)\, 
\frac {\gamma_s^2}{2s!}\,\frac 1{n^{s-2}} + O\big(n^{-(s-1)}\big).
\en
}

\vskip2mm
Using (11.4), one can evaluate the integrals in (11.5) and rewrite them as polynomials
in the cumulants $\gamma_3,\dots,\gamma_{2j+1}$, which in turn may be expressed
polynomially in terms of the moments $\alpha_r = \E X^r$, $r \leq 2j+1$.

\vskip5mm
{\bf Proof.}
The representation (11.3) with $k=2s$ may be written as
$$
\frac{p_n(x)}{\varphi(x)} \, = \, 1 + R_n(x) + 
\frac{\ep_n(x)}{n^{s-1}}\ \frac{1}{\varphi(x)(1+|x|^{2s})}, \qquad
R_n(x) = \sum_{\nu=1}^{2s-2} \frac{q_{\nu}(x)}{n^{\nu/2}},
$$
where $\sup_x |\ep_n(x)| = o(1)$ as $n \rightarrow \infty$.
Since every polynomial $q_\nu$ has degree at most $r = 3(2s-2)$, we necessarily have
$|R_n(x)| \leq \frac{C}{\sqrt{n}}\,(1 + |x|^r)$ up to some constant $C$.
It follows that
$$
\bigg|\frac{p_n(x)}{\varphi(x)} - 1\bigg| \, \leq \, 
\frac{C}{\sqrt{n}}\,(1 + |x|^r) +
\frac{o(n^{-(s-1)})}{\varphi(x)(1+|x|^{2s})} \, \leq \, \delta_n \, \rightarrow 0
$$
as $n \to \infty$ uniformly in $|x|\le M_n(s)$. Using the Lipschitz property 
of the power function near the point 1, we thus obtain that
$$
\bigg(\frac{p_n(x)}{\varphi(x)}\bigg)^\alpha = \,
(1 + R_n(x))^\alpha + \frac{o(n^{-(s-1)})}{\varphi(x)(1+|x|^{2s})},
$$
so that
\be
I(M_n(s)) = \int_{|x|\le M_n(s)} (1 + R_n(x))^\alpha\,\varphi(x)\,dx + o\big(n^{-(s-1)}\big).
\en

Using a Taylor expansion of $x^\alpha$ yields
$$
(1 + R_n(x))^\alpha = 
1 + \sum_{m=1}^{2s-2} 
\frac{\alpha (\alpha - 1) \dots (\alpha - m + 1)}{m!}\,R_n(x)^m +
\frac{C_n(x)}{n^{s-1/2}}\,\big(1 + |x|^{r(2s-1)}\big)
$$
with $\sup_x |C_n(x)| \leq C$ (where $C$ is a constant). Thus, integration in (11.8) 
leads to
$$
I(M_n(s)) = 1 + \sum_{m=1}^{2s-2} 
\frac{\alpha (\alpha - 1) \dots (\alpha - m + 1)}{m!}\,
\int_{|x|\le M_n(s)}R_n(x)^m\,\varphi(x)\,dx + o\big(n^{-(s-1)}\big).
$$
Here the integrals may be extended to the whole real line at the expense 
of an error at most 
$o(\frac{1}{n^{s-1}})$. Indeed, with some constant $C_l$ depending on 
$\l \geq 1$, we have
$$
\int_{|x|>M_n(s)} |x|^l \varphi(x)dx \, \leq \, C_l M_n(s)^{l-1} e^{-M_n(s)^2/2} = 
O\Big(\frac{\log^{\frac{l-1}{2}} n}{n^{s-1}}\Big),
$$
which may be used in the polynomial bound on $R_n$ (together with the factor
$1/\sqrt{n}$). Thus,
$$
I(M_n(s)) = 1 + \sum_{m=1}^{2s-2} 
\frac{\alpha (\alpha - 1) \dots (\alpha - m + 1)}{m!}\,
\int_{-\infty}^\infty R_n(x)^m\,\varphi(x)\,dx + o\big(n^{-(s-1)}\big).
$$
Using a multinomial expansion, we get
$$
R_n(x)^m \ = \sum_{m_1 + \dots + m_{2s-2} = m} \frac{m!}{m_1! \dots m_{2s-2}!} 
n^{-N/2}\,
q_1(x)^{m_1} \dots q_{2s-2}(x)^{m_{2s-2}},
$$
where $N = m_1 + 2m_2 + \dots + (2s-2)\,m_{2s-2}$. That is,
up to a $o(n^{-(s-1)})$--term, one can describe $I(M_n(s)) - 1$ as the sum
\be
\sum \frac{\alpha (\alpha - 1) \dots (\alpha - m + 1)}{m_1! \dots m_{2s-2}!}\, 
n^{-N/2} \,
\int_{-\infty}^\infty q_1(x)^{m_1} \dots q_{2s-2}(x)^{m_{2s-2}}\,\varphi(x)\,dx,
\en
where the summation extends over all integers $m_1, \dots, m_{2s-2} \geq 0$,
not all zero, such that $m = m_1 + m_2 + \dots + m_{2s-2} \leq 2s-2$.

This representation simplifies thanks to the following property of Hermite polynomials:
$$
\int_{-\infty}^\infty H_{\nu_1}(x) \dots H_{\nu_k}(x)\,\varphi(x)\,dx = 0 \qquad
(\nu_1 + \dots + \nu_k \ {\rm is \ odd}).
$$
Hence, it follows from (11.4) that a similar property holds for $q_j$'s as well, so that
the integral in (11.9) is vanishing, as long as $N$ is odd. Restricting ourselves to
the values $N = 2j$, we necessarily have $m_l = 0$ for $l > 2j$, and (11.8) becomes
\be
\sum \frac{\alpha (\alpha - 1) \dots (\alpha - m + 1)}{m_1! \dots m_{2j}!}\, n^{-j} \,
\int_{-\infty}^\infty q_1(x)^{m_1} \dots q_{2j}(x)^{m_{2j}}\,\varphi(x)\,dx,
\en
where the summation extends over all $m_1, \dots, m_{2j} \geq 0$ such that
$m_1 + 2m_2 + \dots + 2j\,m_{2j} = 2j$ and with $m = m_1 + \dots + m_{2j}$.
Finally, we may exclude the case $m_{2j} = 1$, $m_l = 0$ for $l < 2j$, where 
again the above integral is vanishing. As a result, we arrive at the required 
expansion (11.5) with coefficients (11.6). Finally, in the second assertion, 
we necessarily have $b_j=0$ for $j=1,\dots,s-3$ and $b_{s-2} = \gamma_s^2/s!$ 
and we obtain (11.7).
\qed

\vskip5mm
Note that the integral in (11.10) is zero as well, provided that $m=1$ 
(i.e., only one $m_l = 1$). For the index $\alpha = 2$, the factor in front 
of the integral in (11.6) is vanishing unless $m \leq 2$. Hence, we are reduced 
to tuples $m_1,\dots,m_{2j-1}$ such that $m_l = 1$ holds for two different indexes, 
say, $l = \nu_1$ and $l = \nu_2$, and also for tuples where $m_l=2$ holds 
for one $l$ only. Hence, the description of 
the coefficients may be simplified to
$$
b_j = \sum_{\scriptstyle \nu_1,\nu_2>0 \atop \scriptstyle \nu_1+\nu_2=2j}
\int_{-\infty}^\infty q_{\nu_1}(x) q_{\nu_2}(x)\,\varphi(x)\,dx \qquad (\alpha = 2).
$$

Recall that if $T_\alpha(Z_n||Z)$ is finite, then $\E\,e^{cZ_n^2} < \infty$, and hence
$\E\,e^{cX^2} < \infty$ for some $c>0$ (so that $X$ has finite moments of all orders).
In addition, $Z_n$ must have a density in $L^2$. Therefore, all conditions of 
Lemma 11.1 are fulfilled, and in view of the lower bound (11.1), Lemma 11.2 yields:

\vskip5mm
{\bf Proposition 11.3.} {\it For every fixed $s=3,4,\dots$, we have, as $n\to\infty$,
$$
T_\alpha(Z_n||Z) \, \ge \, \frac{1}{\alpha - 1} 
\sum_{j=1}^{s-2}\frac {b_j}{n^j} + O\Big(\frac {1}{n^{s-1}}\Big)
$$
with coefficients given in $(11.6)$. In particular, if $\gamma_j=0$ for 
$j=3,\dots,s-1$ and $\gamma_s\ne 0$, then
\be
T_\alpha(Z_n||Z) \, \ge \, \alpha
\frac{\gamma_s^2}{2s!}\,\frac 1{n^{s-2}} + O\Big(\frac {1}{n^{s-1}}\Big).
\en
}

\vskip2mm
The last lower bound extends to $D_\alpha$ as well
(which is equivalent to $T_\alpha$ when these two distances are small).
Hence we get:

\vskip5mm
{\bf Corollary 11.4.} {\it If, for some integer $K>1$, 
$$
\liminf_{n\to\infty} \, \frac{\log D_\alpha(Z_n||Z)}{\log n}<-K, 
$$
then $\gamma_j=0$ for all $j=3,\dots,K$. In particular,
the random variable $X$ is standard normal, if and only if
$$
\liminf_{n\to\infty}\, \frac{\log D_\alpha(Z_n||Z)}{\log n}=-\infty.
$$
}

\vskip2mm
Combining the lower bound (11.2) with the upper bound (9.2) yields:

\vskip5mm
{\bf Corollary 11.5.}
{\it Let $D_\alpha(X||Z)<\infty$, with $\gamma_j=0,\,j=3,\dots,s-1$, and $\gamma_s\ne 0$
for some $s \geq 3$. Then as $n\to\infty$
$$
\Big(1+O\Big(\frac 1n\Big)\Big)\,
\frac{\gamma_s^2}{2s!}\,\frac{1}{n^{s-2}}\le D_\alpha(Z_n||Z)\le n D_\alpha(X||Z).
$$
}

\vskip10mm
\section{{\bf Necessity Part in Theorem 1.2 ($d=1$)}}
\setcounter{equation}{0}

\vskip2mm
\noindent
Again, let $X,X_1,X_2,\dots$ denote i.i.d. random variables with characteristic
function $f(t) = \E\,e^{itX}$, and let $Z_n = (X_1 + \dots + X_n)/\sqrt n$.
The necessity part in Theorem 1.2 does not require any moment assumptions 
on the mean and variance.
As a preliminary step, the next proposition provides a subgaussian bound 
on the Laplace transform $f(iy) = \E\,e^{-yX}$ subject to the sublinear growth 
of $D_\alpha(Z_n||Z)$.
Recall that $\alpha > 1$ is fixed, and we denote its conjugate value by 
$\beta = \alpha/(\alpha - 1)$.

\vskip5mm
{\bf Lemma 12.1.}
{\it If\, $\liminf_{n\to\infty} \big[\frac{1}{n}\,D_\alpha(Z_n||Z)\big] = 0$, then 
\be
f(iy) \le e^{\beta y^2/2}, \qquad y\in\mathbb R. 
\en
}

\vskip2mm
{\bf Proof.}
Indeed, by Proposition 4.2, applied to $Z_n$ in place of $X$, for all $y \in \R$,
$$
f(iy/\sqrt n)^n \, \le \, 
\big(1 + (\alpha - 1)T_\alpha(Z_n||Z)\big)^{1/\alpha}\, e^{\beta y^2/2},
$$
and after a change of the variable we get
$$
f(iy) \le \exp\Big\{\frac{1}{\alpha n}\,\log
\big(1 + (\alpha - 1)T_\alpha(Z_n||Z)\big)\Big\}\, e^{\beta y^2/2}.
$$
But $\liminf_{n\to\infty} \big[\frac{1}{n}\,D_\alpha(Z_n||Z)\big] = 0$, if and only if
$\liminf_{n\to\infty} \big[\frac{1}{n}\,\log(1 + (\alpha - 1)T_\alpha(Z_n||Z))\big] = 0$.
Hence, we arrive at the required conclusion by letting $n \rightarrow \infty$ along
a suitable subsequence.
\qed

\vskip5mm
In other words, if $f(iy_0)>e^{\beta y_0^2}/2$ holds for some $y_0\in\mathbb R$, then 
$D_\alpha(Z_n||Z)\ge cn$ holds for some positive constant $c$. 
Thus, in this case $D_\alpha(Z_n||Z)$ has a maximal growth rate, in view of
the sublinear upper bound (9.2).

The assumption of Lemma 12.1 is fulfilled, when $D_\alpha(Z_n||Z) \rightarrow 0$, which
provides a slightly weakened variant of the necessary condition (1.4) in Theorem 1.2
for dimension $d=1$ (replacing the strict inequality with a non-strict inequality). 
To arrive at a more precise condition, we have to add another preliminary step.

\vskip5mm
{\bf Lemma 12.2.}
{\it If\, $\lim_{n\to\infty} D_\alpha(Z_n||Z)\big] = 0$, then, for any integer 
$k \geq \alpha/2$,
\be
\lim_{n\to\infty} \, 
\int_{-\infty}^\infty f\big(iy/\sqrt{kn}\big)^{2kn}\, e^{-\beta y^2}\,dy \, = \,
\sqrt{\pi (\alpha - 1)}.
\en
}

\vskip2mm
{\bf Proof.} The characteristic function of $Z_n$ is given by $f_n(t) = f(t/\sqrt{n})^n$.
Hence, the integral in (12.2) is just
\bee
\int_{-\infty}^\infty \big(\E\, e^{-yZ_{nk}}\big)^2\, e^{-\beta y^2}\,dy 
 & = &
\int_{-\infty}^\infty \E\, e^{-y\,(Z_{nk} + Z_{nk}')}\, e^{-\beta y^2}\,dy \\
 & = &
\int_{-\infty}^\infty \E\, e^{- \sqrt{2}\,yZ_{2nk}}\, e^{-\beta y^2}\,dy \, = \,
\sqrt{\frac{\pi}{\beta}}\, \E\,e^{\frac{1}{2\beta} Z_{2nk}^2},
\ene
where by $Z_{nk}'$ we denoted an independent copy of $Z_{nk}$.
On the other hand, since $Z_{nk}$ is a normalized sum of $k$ independent
copies of $Z_n$, we may apply Proposition 5.1 with $X$ replaced by $Z_n$ and
with $n$ replaced by $2k$. In this case, inequality (5.2) tells us that, whenever 
$2k \geq \alpha$, we have
$$
\big|\E\,e^{\frac{1}{2\beta} Z_{2nk}^2} - \E\,e^{\frac{1}{2\beta} Z^2}\big|  \, \leq \,
 c_{2k}\, \Big(\big(1 + \chi_\alpha(Z_n||Z)^{1/\alpha}\big)^k - 1\Big), \qquad
Z \sim N(0,1).
$$
Since, by the assumption, $\chi_\alpha(Z_n,Z) \rightarrow 0$ as $n \rightarrow \infty$,
the limit in (12.2) is equal to 
$\sqrt{\frac{\pi}{\beta}}\, \E\,e^{\frac{1}{2\beta} Z^2}$.
\qed

\vskip2mm
{\bf Proof of the neccesity part in Theorem 1.2 for $d=1$.} 
Let $D_\alpha(Z_n||Z) \rightarrow 0$ as $n \rightarrow \infty$.
Given a fixed number $\delta > 0$, let us decompose
\begin{eqnarray}
\hskip-15mm
\int_{-\infty}^\infty f\big(iy/\sqrt{nk}\big)^{2nk}e^{-\beta y^2}\,dy
 & = &
I_1+I_2 \nonumber \\
 & = &
\bigg(\int_{|y|\le \delta\sqrt{nk}}+\int_{|y|>\delta\sqrt{nk}}\bigg)
f\big(iy/\sqrt{nk}\big)^{2nk} e^{-\beta y^2}\,dy.
\end{eqnarray}
The characteristic function $f$ is entire, and $f(0)=1$, hence it is non-vanishing in 
some disc $|t| < R$ on the complex plane. Define $g(t)=\log f(t)$ for $|t|<R$, choosing 
the branch of the logarithm according to  the condition $\log f(0)=0$. The function $g$ 
is analytic in the same disc and admits a power series representation
$$
g(t)=-\frac 12 t^2+\sum_{k=3}^{\infty}b_k t^k.
$$
Clearly, for a suitable value $r \in (0,R)$ and some constant $C$, we have 
$\sum_{k=3}^{\infty} |b_k t^k| \leq C |t|^3$ in the disc $|t| \leq r$,
so that
$$
f\big(iy/\sqrt{nk}\big)^{2nk} = \exp\{y^2 + \theta y^3/\sqrt n\}\quad
\text{ for real} \ \ y\in[-r\sqrt{nk},r\sqrt{nk}\,],
$$
where $\theta$ is a quantity such that $|\theta|\le C$ and $k \geq \alpha/2$ is a fixed
integer. Assuming that $\delta \leq \min\{r,1/(2C)\}$, this relation allows us 
to rewrite the integral $I_1$ as
$$
I_1 = \int_{|y|\le \delta\sqrt{nk}} e^{-(\beta - 1) y^2 + \theta y^3/\sqrt n}\,dy.
$$
Here the term $\theta y^3/\sqrt n$ in the above exponent may be removed at the expense 
of an error of order $O(\frac{1}{\sqrt{n}})$. This is justified by the bounds
\bee
\int_{|y|\le \delta\sqrt{nk}} 
\big|e^{-(\beta - 1) y^2 + \theta y^3/\sqrt n} - e^{-(\beta - 1) y^2}\big|\,dy
 & \leq &
\int_{|y|\le \delta\sqrt{nk}} \frac{Cy^3}{\sqrt n}\, e^{-(\beta - 1) y^2 + 
C |y|^3/\sqrt n}\, dy \\
 & \leq &
\frac{C}{\sqrt n}
\int_{|y|\le \delta\sqrt{nk}} y^3\, e^{-(\beta - 1) y^2/2}\, dy \, = \, 
O\Big(\frac{1}{\sqrt{n}}\Big).
\ene
Hence
$$
I_1 = \int_{|y|\le \delta\sqrt{nk}} e^{-(\beta - 1) y^2} dy + O\Big(\frac{1}{\sqrt{n}}\Big)
= \sqrt{\pi (\alpha - 1)} + O\Big(\frac{1}{\sqrt{n}}\Big),\qquad n\to\infty.
$$
Applying this result in (12.3), the equality (12.2) implies that
$I_2 \rightarrow 0$, or equivalently
\be
\int_{|u|>\delta} \big(f(iu)\,e^{-\beta u^2/2}\big)^{2nk}\,du = 
o\Big(\frac 1{{\sqrt n}}\Big)\quad
\text{as}\ n\to\infty,
\en
which holds for any sufficiently small $\delta > 0$, and since the integrand is 
non-negative, for any smaller fixed $\delta>0$ as well.

Now, the function $\psi(u) = f(iu) \,e^{-\beta u^2/2}$ is analytic and satisfies 
$0 < \psi(u) \leq 1$ on the real line, cf. (12.1). In order to show that $\psi(u) < 1$ 
for all $u \neq 0$, suppose for a moment that $\psi(u_0) = 1$ for some $u_0 > 0$. 
Obviously $u_0$ has to be local maximum point, which implies $\psi'(u_0) = 0$. 
Hence the power series representation at this point, that is
$$
\psi(u) - 1 = c_l (u - u_0)^l + \sum_{j=l+1}^\infty c_j (u - u_0)^j
$$
starts with a non-zero term $c_l \neq 0$ for some $l \geq 2$. Since $\psi(u) - 1 \leq 0$ 
for all $u \in \R$, all coefficients are real numbers, and moreover, $l = 2m$ is even 
($m \geq 1$) and 
$c_l < 0$. Hence, in some neighborhood $|u - u_0| \leq r_0 < u_0$ and for some constants 
$c_1,c_0>0$, we have $\psi(u) \geq 1 - c_1 (u - u_0)^{2m} \geq e^{-c_0  (u - u_0)^{2m}}$.
Now choosing $\delta = u_0 - r_0$, this neighborhood is contained in $(\delta,\infty)$,
and with some constant $c>0$ we get
\bee
\int_{|u|>\delta} \big(f(iu)\,e^{-\beta u^2/2}\big)^{2nk}\,du 
 & \geq &
\int_{|u-u_0| < \delta} \psi(u)^{2nk}\,du \\
 & \geq &
\int_{|u-u_0| < \delta} \exp\big\{-2nk\cdot c_0  (u - u_0)^{2m}\big\}\,du \\
 & = &
2 \int_0^\delta \exp\big\{-2nk \cdot c_0  x^{2m}\big\}\,dx \ \geq \ \frac{c}{n^{1/(2m)}},
\ene
which contradicts to the asymptotic relation (12.4). The case $u_0 < 0$ is similar, and 
thus we necessarily arrive at $\psi(u) < 1$ for all real $u \neq 0$.
\qed

\vskip10mm
\section{{\bf Pointwise Upper Bounds for Convolutions of Densities}}
\setcounter{equation}{0}

\vskip2mm
\noindent
Before turning to the sufficiency part in Theorem 1.2,  we shall derive several upper bounds 
for the densities $p_n$ of the normalized sums $Z_n$. In general, bounds for the density 
$p(x)$ of $X$ at individual points $x$ cannot be deduced from $D_\alpha(X||Z) < \infty$.
However, this is possible after several convolutions of $p$ with itself. The following 
observation holds without assuming that $X$ has mean zero and variance one. Put
\bee
\psi(u) 
 & = &
f(iu) e^{-\beta u^2/2} \\
 & = &
\E\,e^{-uX}\,e^{-\beta u^2/2} \, = \, 
e^{-\beta u^2/2} \int_{-\infty}^\infty e^{-ux}\,p(x)\,dx, \qquad u \in \R,
\ene
where $f$ is the characteristic function of $X$ and $\beta = \frac{\alpha}{\alpha - 1}$.
As usual, $Z$ denotes a standard normal random variable.

\vskip5mm
{\bf Proposition 13.1.} {\it Given a random variable $X$ such that 
$T_\alpha = T_\alpha(X||Z) < \infty$, 
we have, for all $x \in \R$ and $n \geq n_\beta = \max(\beta,2)$,
\be
p_n(x) \, \le \, \frac{A_\alpha \sqrt{n}}{(2\pi)^{1/2}}\, e^{-x^2/(2\beta)}\,
\psi\Big(-\frac{x}{\beta\sqrt n}\Big)^{n-n_\beta},
\en
where $A_\alpha = \big(1 + (\alpha - 1)T_\alpha\big)^{\frac{1}{\alpha - 1}}$ in case 
$1 < \alpha \leq 2$ and
$A_\alpha = \big(1 + (\alpha - 1)T_\alpha\big)^{\frac{2}{\alpha}}$ in case $\alpha > 2$.
}

\vskip5mm
In particular, under the condition (1.1), that is, when $\psi \leq 1$, 
we arrive at the following subgaussian pointwise bound
$$
p_n(x) \, \le \, \frac{A_\alpha \sqrt{n}}{(2\pi)^{1/2}}\, e^{-x^2/(2\beta)},
$$
which may be effective in the region $|x| >\!> \sqrt{\log n}$. It can be sharpened further
for larger values of $|x|$ by virtue of Proposition 4.3. Combined with (13.1), 
it immediately provides an exponential pointwise bound (with respect to $n$).

\vskip5mm
{\bf Corollary 13.2.} {\it If $T_\alpha(X||Z) < \infty$, there exist
constants $x_0 > 0$ and $\delta \in (0,1)$ depending on the density $p$ only, such that,
for all $n$ large enough,
\be
p_n(x) \, \le \, \delta^n  e^{-x^2/(2\beta)}\,\psi\Big(-\frac{x}{\beta\sqrt n}\Big)^{n/2}, 
\quad 
whenever \ \ |x| \geq x_0 \sqrt{n}.
\en
}

\vskip2mm
Here the last $\psi$-term is (13.2) will become crucial for bounding $T_\alpha(Z_n||Z)$.

\vskip5mm
{\bf Proof of Proposition 13.1.}
Since $\E\, e^{yX}<\infty$ for all $y\in\R$, the characteristic function 
$f_n(t) = \E\,e^{itZ_n} = f(t/\sqrt{n})^n$ is extended as an entire function 
to the complex plane. Since $p$ belongs to $L^\alpha(\R,dx)$, an application 
of the Hausdorff-Young inequality
implies that $f_n$ is integrable whenever $n \geq \max(\beta,2)$.
In this case $Z_n$ has a continuous density given by the Fourier inversion formula 
$$
p_n(x) \, = \, \frac{1}{2\pi}\int_{-\infty}^\infty e^{-itx} f(t/\sqrt n)^n\,dt \, = \,
\frac{1}{2\pi}\, \lim_{T \rightarrow \infty} \int_{-T}^T e^{-itx} f(t/\sqrt n)^n\,dt.
$$
Moreover, since the family $\{e^{hx} p_n(x)\}_{0 \leq h \leq y}$ is compact in 
$L^1(\R)$,
$f_n(t)$ tends to zero at infinity uniformly in every strip $|\Im\, t|\le y <\infty$
(by the Riemann-Lebesgue lemma). Applying Cauchy's theorem to rectangle contour 
$[-T,T] \cup [T,T+iy] \cup [T+iy,-T+iy] \cup [-T+iy,-T]$, the inversion formula 
may therefore be written as
\be
p_n(x) = e^{yx} \frac{1}{2\pi}\int_{-\infty}^\infty e^{-itx}f((t + iy)/\sqrt n)^n\,dt
\en
for any fixed $y>0$. Without loss of generality, let $x<0$. 

{\bf Case} $\alpha > 2$, $n \geq 2$. Using
$|f(t+iy)|\le f(iy)$ ($t,y\in\mathbb R$) and changing variable in (13.3), we get
\be
p_n(x) \, \le \,
e^{yx}f(iy/\sqrt n)^{n-2}\, \frac{\sqrt n}{2\pi} 
\int_{-\infty}^\infty |f(t + iy/\sqrt n)|^2\,dt.
\en
The function $t \rightarrow f(t + iy/\sqrt n) = \E\,e^{itX - yX/\sqrt{n}}$ is the Fourier 
transform of $g(u) = e^{-yu/\sqrt n}\, p(u)$. Hence, by Parseval's identity,
$$
\frac{1}{2\pi} \int_{-\infty}^\infty |f(t + iy/\sqrt n)|^2\,dt =
\int_{-\infty}^\infty e^{-2yu/\sqrt n}\, p(u)^2\,du.
$$
To estimate the latter integral, factorize the integrand as
$\big(e^{-2yu/\sqrt n}\varphi(u)^{2/\beta}\big)\, \frac{p(u)^2}{\varphi(u)^{2/\beta}}$
and apply H\"older's inequality with exponents $r = \frac{\alpha}{\alpha - 2}$,
$r^* = \frac{\alpha}{2}$. Thus, up to the factor 
$(1 + (\alpha - 1)T_\alpha)^{2/\alpha}$,
this integral can be estimated from above by
$$
\bigg(\int_{-\infty}^\infty e^{-2ryu/\sqrt n}\,
\varphi(u)^{2r/\beta}\,du\bigg)^{1/r} \, = \,
\frac{1}{\sqrt{2\pi}}\,\Big(\frac{\alpha - 2}{2\alpha - 2}\Big)^{1/2r}\,e^{\beta y^2/n} 
\, \leq \,
\frac{1}{\sqrt{2\pi}}\,e^{\beta y^2/n}.
$$
This gives
$$
\frac{1}{2\pi} \int_{-\infty}^\infty |f(t + iy/\sqrt n)|^2\,dt \, \le \, 
\frac{1}{\sqrt{2\pi}}\, \big(1 + (\alpha - 1)T_\alpha\big)^{2/\alpha}\,e^{\beta y^2/n},
$$
and (13.4) results in the upper bound
\bee
p_n(x)
 & \le &
\sqrt{\frac{n}{2\pi}}\ \big(1 + (\alpha - 1)T_\alpha\big)^{2/\alpha}\, 
e^{yx+\beta y^2/n}\,f(iy/\sqrt n)^{n-2} \\
 & = &
\sqrt{\frac{n}{2\pi}}\ \big(1 + (\alpha - 1)T_\alpha\big)^{2/\alpha}\, 
e^{yx+\beta y^2/2}\, \psi(y/\sqrt n)^{n-2}.
\ene
Choosing here $y=-x/\beta$, we arrive at (13.1).

\vskip2mm
{\bf Case} $1 < \alpha \leq 2$, $n \geq \beta$. Again using
$|f(t+iy)|\le f(iy)$ ($t,y\in\mathbb R$) and changing variable, we obtain from 
(13.3) that
\be
p_n(x) \, \le \,
e^{yx}f(iy/\sqrt n)^{n-\beta}\, 
\frac{\sqrt n}{2\pi} \int_{-\infty}^\infty |f(t + iy/\sqrt n)|^\beta\,dt.
\en
Now, since $\beta \geq 2$, we are allowed to apply the classical Hausdorff-Young 
inequality
$$
\bigg(\frac{1}{2\pi} \int_{-\infty}^\infty 
|f(t + iy/\sqrt n)|^\beta\,dt\bigg)^{1/\beta} \leq
\|g\|_\alpha = 
\bigg(\int_{-\infty}^\infty e^{-\alpha yu/\sqrt n}\, p(u)^\alpha\,du\bigg)^{1/\alpha}.
$$
To estimate the latter integral, factorize its integrand as 
$\big(e^{-\alpha yu/\sqrt n}\,
\varphi(u)^{\alpha - 1}\big)\, \frac{p(u)^\alpha}{\varphi(u)^{\alpha - 1}}$ 
and use the inequality
$e^{-\alpha yu/\sqrt n}\,\varphi(u)^{\alpha - 1} \le 
(2\pi)^{-\frac{\alpha - 1}{2}}\,e^{\frac{\alpha \beta y^2}{2n}}$. This gives
$$
\frac{1}{2\pi} \int_{-\infty}^\infty |f(t + iy/\sqrt n)|^\beta\,dt \leq 
\Big((2\pi)^{-\frac{\alpha - 1}{2}}\,e^{\frac{\alpha \beta y^2}{2n}}\,
(1 + (\alpha - 1)T_\alpha)\Big)^{\frac{\beta}{\alpha}} = 
\frac{1}{\sqrt{2\pi}}\,e^{\frac{\beta^2 y^2}{2n}}\,
(1 + (\alpha - 1)T_\alpha)^{\frac{1}{\alpha - 1}}.
$$
Hence, (13.5) results in the upper bound
\bee
p_n(x)
 & \le &
\sqrt{\frac{n}{2\pi}}\,\big(1 + (\alpha - 1)T_\alpha\big)^{\frac{1}{\alpha - 1}}\, 
e^{yx+\beta^2 y^2/2n}\, f(iy/\sqrt n)^{n-\beta}  \nonumber \\
 & = &
\sqrt{\frac{n}{2\pi}}\,\big(1 + (\alpha - 1)T_\alpha\big)^{\frac{1}{\alpha - 1}}\, 
e^{yx+\beta y^2/2}\, \psi(y/\sqrt n)^{n-\beta}.
\ene
Again choosing $y=-x/\beta$, we arrive at (13.1).
\qed

\vskip10mm
\section{{\bf Sufficiency Part in Theorem 1.2 ($d=1$).}}
\setcounter{equation}{0}

\vskip2mm
\noindent
Let $X,X_1,X_2,\dots$ be i.i.d. random variables such that $\E X = 0$, $\E X^2 = 1$, 
with characteristic function $f(t) = \E\,e^{itX}$. As before, put 
$\psi(u)=f(iu)\,e^{-\beta u^2/2}$, $\beta = \frac{\alpha}{\alpha - 1}$, and let 
$Z \sim N(0,1)$.
Assuming that the condition (1.4) is fulfilled, i.e., $\psi(u) < 1$ for all real 
$u \neq 0$, here it will be shown that the normalized sums
$$
Z_n = \frac{X_1 + \dots + X_n}{\sqrt n}
$$
do satisfy $T_\alpha(Z_n||Z) \rightarrow 0$ as $n \rightarrow \infty$, as long as 
$T_\alpha(Z_{n_0}||Z) < \infty$ for some $n_0$. We also derive an asymptotic 
expansion for this distance which is rather similar to (1.2) in case $\alpha = 2$. 
For simplicity, let us assume that $n_0 = 1$,  so that $X$ has density $p$ with 
$T_\alpha(X||Z) < \infty$ (the general case $n_0 \geq 1$ 
is rather similar and needs only minor modifications). In particular, all $Z_n$ have 
densities $p_n$ which are continuous and bounded for all $n$ large enough.

In Section 11, we considered integrals of the form
$$
I_0 = \int_{|x|\le M_n} \frac{p_n(x)^\alpha}{\varphi(x)^{\alpha - 1}}\,dx \qquad 
{\rm with} \ \ 
M_n = \sqrt{2(l-1)\log n} \ \ (l=3,4,\dots)
$$
According to Proposition 11.1 with $k=2l$ and Lemma 11.2, these integrals admit 
an asymptotic expansion
\be
I_0 = 1 + \sum_{j=1}^{l-1} \frac{b_j}{n^j} + o\big(n^{-(l-1)}\big),
\en
which may be simplified in terms of the cumulants of $X$ as
\be
I_0 = 1 + \alpha (\alpha-1)\, 
\frac {\gamma_s^2}{2s!}\,\frac 1{n^{s-2}} + O\big(n^{-(s-1)}\big)
\en
when $\gamma_j=0$ for $j=3,\dots,l-1$. Hence, for the proof of Theorem 1.2
(in dimension one), it remains to bound the integral of $p_n^\alpha/\varphi^{\alpha - 1}$ 
over the complementary region $|x| > M_n$ by a polynomially small quantity 
(with respect to $n$). More precisely, it will be sufficient to show that, 
for any large enough $l\ge 3$ and some constant $\kappa>0$,
\be
\int_{|x|>M_n}\frac{p_n(x)^\alpha}{\varphi(x)^{\alpha - 1}}\,dx \, = \, 
O\Big(\frac 1{n^{\kappa l}}\Big),
\qquad n\to\infty.     
\en

To this aim, we need to properly estimate $p_n(x)$ for $|x|>M_n$, which can be done
based on the pointwise bounds of the previous section. For definiteness, let us
consider the half-axis $x<-M_n$, which we split into three intervals reflecting 
the possible different behavior of these densities. Namely, define
$$
I_1  =
\int_{-\infty}^{-x_0\sqrt n} \frac{p_n(x)^\alpha}{\varphi(x)^{\alpha - 1}}\,dx, \qquad
I_2 = \int_{-x_0\sqrt n}^{-x_1\sqrt n} \frac{p_n(x)^\alpha}{\varphi(x)^{\alpha - 1}}\,dx, 
\qquad
I_3 = \int_{-x_1\sqrt n}^{-M_n}\,  \frac{p_n(x)^\alpha}{\varphi(x)^{\alpha - 1}}\,dx
$$ 
with parameters $0 < x_1 < x_0$ and assuming that $M_n < x_1 \sqrt n$
(otherwise, $I_3 = 0$).

Applying inequality (13.2), we obtain that, for all $n$ large enough and with some 
$\delta \in (0,1)$ and $x_0 > 0$,
\bee
I_1 
 & \leq &
(2\pi)^{\frac{\alpha - 1}{2}}\,\delta^{\alpha n} 
\int_{-\infty}^{-x_0\sqrt n} \psi\Big(-\frac{x}{\beta\sqrt n}\Big)^{\alpha n/2}\, dx \\
 & \leq &
(2\pi)^{\frac{\alpha - 1}{2}}\,\delta^{\alpha n} \beta \sqrt{n}
\int_{-\infty}^{\infty} \psi(u)^m\, du, \qquad m \leq \frac{\alpha n}{2},
\ene
where on the last step we used $\psi \leq 1$. By Corollary 5.2, cf. (5.3), 
the last integral is convergent whenever $m \geq \alpha$. One may take, for example, 
$m = [\alpha] + 1$, which ensures the condition $m \leq \frac{\alpha n}{2}$ 
for all sufficiently large $n$. Hence
$$
I_1 \, \le \, C\delta_1^n \qquad (n \geq n_1)
$$
with some constants $C>0$, $x_0>0$ and $\delta < \delta_1 < 1$,
depending on the density $p$ only.

To estimate the integral $I_2$ (with any fixed number $0< x_1 < x_0$), we employ 
Proposition 13.1. By the condition (1.4), the function $\psi$ is bounded away 
from 1 on any compact interval in $(-\infty,0)$, so, 
$\delta_2 = \max_{-x_0 \leq u \leq -x_1} \psi(u) < 1$. 
Hence, by inequality (13.1),
\bee
I_2 
 & \le & 
A_\alpha n^{\alpha/2}
\int_{-x_0\sqrt n}^{-x_1\sqrt n} \psi\Big(-\frac{x}{\beta\sqrt n}\Big)^{n-n_\beta}\,dx \\
 & = &
A_\alpha^\alpha\, \beta\, n^{(\alpha + 1)/2} 
\int_{-x_0/2}^{-x_1/2} \psi(u)^{n - n_\beta}\,du
 \ \le \ 
A_\alpha^\alpha\, \beta\, n^{(\alpha + 1)/2}\, (x_0-x_1)\,\delta_2^{n - n_\beta}
\ene
which again decays exponentially fast like $I_1$.

It remains to properly estimate the integral $I_3$ with some (prescribed) $x_1>0$.
In order to estimate $p_n(x)$ in $[-x_1\sqrt n,-M_n]$, we use the bound (13.1) once more.
As discussed in Section 12, the function $h(u) = \log f(iu)$ is analytic in some 
disc $|u| \leq r$, and since $h(0) = 0$, $h'(0) = 1/2$, we have 
$h(u) \sim \frac{1}{2}\,u^2$
near zero. Hence $|h(u)| \leq \frac{1+\beta}{4}\,|u|^2$ throughout this disc, when 
$r$ is sufficiently small, implying $|f(iu)| \leq e^{(1+\beta) |u|^2/4}$. Hence for 
$u$ real, $|u| \leq r$, 
we have $\psi(u) \leq e^{-(\beta - 1) |u|^2/4}$, which implies
$$
\psi\Big(-\frac{x}{\beta\sqrt n}\Big)^{n-n_\beta} \leq 
\psi\Big(-\frac{x}{\beta\sqrt n}\Big)^{n/2} \leq 
\exp\Big\{-\frac{\beta - 1}{4}\, \frac{x^2}{2\beta^2}\Big\} = 
e^{-x^2/(8 \alpha \beta)} 
$$
for all $n \geq 2\max(\beta,2)$ and $-\beta r\sqrt n < x < 0$. Therefore, by (13.1), 
in this interval
$$
\frac{p_n(x)^\alpha}{\varphi(x)^{\alpha - 1}} \le 
A_\alpha^\alpha\, n^{\alpha/2}\, e^{-x^2/(8 \beta)},
$$
which results with $x_1 = \beta r$ in
\bee
I_3 
 & \le &
A_\alpha^\alpha\, n^{\alpha/2} \int_{-x_1\sqrt n}^{-M_n} e^{-x^2/(8 \beta)}\,dx \\
 & \leq &
\sqrt{2\pi \beta}\, A_\alpha^\alpha\, n^{\alpha/2}\, e^{-M_n^2/(8 \beta)} \, = \, 
\sqrt{2\pi \beta}\, A_\alpha^\alpha\, n^{-(\frac{l - 1}{4 \beta} - \frac{\alpha}{2})},
\ene
where we used a well-known inequality 
$\int_M^\infty \varphi(x)\,dx \leq \frac{1}{2}\,e^{-M^2/2}$ ($M > 0$).

Collecting these bounds, we obtain that $I_1 + I_2 + I_3 = o(n^{-l/8 \beta})$
for a sufficiently large $l$. A similar relation holds true for integrals over 
the half-axis $x > M_n$, which proves (14.3). 

Since $T_\alpha(Z_n||Z) = \frac{1}{\alpha - 1}\,(I_0 + I_1 + I_2 + I_3 - 1)$, 
and using the expansions (14.1)-(14.2), we conclude that, for any $s = 3,4,\dots$,
\be
T_\alpha(Z_n||Z) = \frac{1}{\alpha - 1} \sum_{j=1}^{s-2} \frac{b_j}{n^j} + 
O\big(n^{-(s-1)}\big)
\en
with coefficients $b_j$ described in (11.6). Moreover, in terms of the cumulants 
of $X$, (14.4) simplifies to
\be
T_\alpha(Z_n||Z) = \alpha\, 
\frac {\gamma_s^2}{2s!}\,\frac 1{n^{s-2}} + O\big(n^{-(s-1)}\big) \quad
{\rm in \ case} \ \ \gamma_j=0 \ {\rm for} \ j=3,\dots,s-1.
\en
Since $D_\alpha$ and $T_\alpha$ are equivalent (when these quantities are small), 
the last relation holds true for the R\'enyi distance $D_\alpha(Z_n||Z)$ as well. 
Thus, Theorem 1.2 is proved in dimension one.
\qed

\vskip10mm
\section{{\bf Non-uniform Local Limit Theorem}}
\setcounter{equation}{0}

\vskip2mm
\noindent
Here we prove Theorem 1.3 in dimension one, still keeping the basic assumptions 
$\E X = 0$, $\E X^2 = 1$. We shall state it in a more precise form, by using 
the cumulants $\gamma_k$ of $X$. We remind that $\beta = \frac{\alpha}{\alpha - 1}$ 
($\alpha > 1$).

\vskip5mm
{\bf Theorem 15.1.} {\it Suppose that $D_\alpha(Z_n||Z)$ is finite for some $n=n_0$, 
and assume that condition $(1.4)$ holds. If $\gamma_3 = \dots = \gamma_{s-1} = 0$ 
for some $s \geq 3$, then
\be
\sup_{x \in \R} \, \frac{|p_n(x) - \varphi(x)|}{\varphi(x)^{1/\beta}} \, = \, 
\frac{a_s\,|\gamma_s|}{s!}\, n^{-\frac{s-2}{2}} + O\big(n^{-\frac{s-1}{2}}\big),
\en
where 
$$
a_s = \sup_{x \in \R}\, \big[\varphi(x)^{1/\alpha}\, |H_s(x)|\big].
$$
}

\vskip2mm
In case $s=3$ we thus obtain the inequality (1.5), and if $\E X^3 = 0$ (and hence
$\gamma_3 = 0$), one may turn to the next moment of order $s = 4$, which yields
the rate $1/n$ in (15.1). As for the cumulant coefficient, let us recall that
$\gamma_s = \E H_s(X) = \E X^s - \E Z^s$ (cf. Proposition 8.1).

To compare these results with Proposition 11.1, note that, assuming the existence of
moments of order $s$, and that $Z_n$ has a bounded continuous density $p_n$ for 
large $n$, the Edgeworth expansion (11.3) allows to derive a weaker statement, 
such as
$$
\sup_{x \in \R} \, (1 + |x|^s)\,|p_n(x) - \varphi(x)| \, = \, 
\frac{a_s'\,|\gamma_s|}{s!}\, n^{-\frac{s-2}{2}} + o\big(n^{-\frac{s-2}{2}}\big),
$$
where $a_s' = \sup_{x \in \R}\, (1 + |x|^s)\, |H_s(x)|\,\varphi(x)$
(still assuming that the moments of $X$ of orders less than $s$ are the same 
as for the standard normal law).

Note in addition that the condition (1.4) is almost necessary for the conclusion 
such as (15.1) and even for a weaker one. Indeed, suppose that
\be
\liminf_{n \rightarrow \infty} \, 
\sup_{x \in \R} \, \frac{p_n(x) - \varphi(x)}{\varphi(x)^{1/\beta}} \, < \, \infty,
\en
so that
$$
p_n(x) \leq \varphi(x) + C_n\, \varphi(x)^{1/\beta}, \qquad 
\liminf_{n \rightarrow \infty} \, C_n < \infty.
$$
Multiplying this inequality by $e^{tx}$ and integrating, we get
$$
\big(\E\,e^{tX/\sqrt{n}}\big)^n = \E\,e^{tZ_n} \leq e^{t^2/2} + BC_n\,e^{\beta t^2/2},
\qquad
B = (2\pi)^{(1-1/\beta)/2} \sqrt{\beta}.
$$
Now substitute $t$ with $t\sqrt{n}$ and raise the above inequality 
to the power $1/n$. Letting $n \rightarrow \infty$ along a suitable subsequence, 
we arrive in the limit at
$$
\E\,e^{tX}\, \leq \, e^{\beta t^2/2}, \qquad t \in \R.
$$
Thus, this subgaussian property is indeed implied by the local limit theorem 
in the form (15.2).

\vskip5mm
{\bf Proof of Theorem 15.1.} Here in contrast with the proof of Theorem 1.2, 
we need to consider a decomposition into a smaller number of zones. 
For simplicity, let $n_0=1$, and as before, define 
$$
M_n = \sqrt{2(l-1)\log n}
$$ 
with parameter $l \geq s-1$, assuming that is sufficiently large. 
Then (11.3) yields the desired equality (15.1), provided that the supremum 
on the left is taken over the interval $|x| \leq M_n$. Hence, it will be 
sufficient to bound the two suprema
$$
J_1  = \sup_{|x| \geq x_1\sqrt n} \, \frac{p_n(x)}{\varphi(x)^{1/\beta}}, \qquad
J_2 = \sup_{M_n \leq |x| \leq x_1\sqrt n}\,  \frac{p_n(x)}{\varphi(x)^{1/\beta}}
$$ 
by polynomially small quantities  (with respect to $n$) 
with some $x_1 > 0$ and assuming that $M_n < x_1 \sqrt n$
(otherwise, $J_2 = 0$).

To this aim, we again invoke the bounds of Proposition 13.1 and Corollary 13.2.
The assumption (1.4) means that the function $\psi(u) = \E\,e^{-uX}\,e^{-\beta u^2/2}$
satisfies $\psi(u)<1$ for all $u \neq 0$. Hence, the bound (13.2) yields, for all $n$ 
large enough,
$$
\frac{p_n(x)}{\varphi(x)^{1/\beta}} \leq \delta^n, \qquad |x| \geq x_0 \sqrt{n},
$$
which is valid with some $\delta \in (0,1)$ and $x_0 > 0$. Moreover, since
$\delta_2 = \max_{x_1 \leq |u| \leq x_0} \psi(u) < 1$ for any $x_1 \in (0,x_0)$, 
the bound (13.1) yields
$$
\frac{p_n(x)}{\varphi(x)^{1/\beta}} \, \le \, A_\alpha \sqrt{n}\, 
\delta_2^{n-n_\beta} \, \leq \, \delta_1^n \qquad (n \geq n_\beta = \max(\beta,2))
$$
with some $\delta < \delta_2 < 1$. Both estimates imply $J_1 = O(\delta_1^n)$
as $n \rightarrow \infty$ for any $x_1 > 0$.

Moreover, as shown in the proof of the sufficiency part of Theorem 1.2, 
we have for some $x_1 > 0$,
$$
\frac{p_n(x)}{\varphi(x)^{1/\beta}} \, \le \, 
A_\alpha \sqrt{n}\,e^{-x^2/(8\alpha \beta)},
\qquad |x| \leq x_1 \sqrt{n}.
$$
This gives 
$$
J_2 \leq A_\alpha \sqrt{n}\,e^{-M_n^2/(8\alpha \beta)} = 
A_\alpha\, n^{-(\frac{l - 1}{8 \alpha \beta} - \frac{1}{2})} \leq 
A_\alpha n^{-\kappa},
$$
where the last inequality holds for any prescribed value of $\kappa>0$ 
by a suitable choice of $l$.
\qed

\vskip10mm
\section{{\bf The Multidimensional Case}}
\setcounter{equation}{0}

\vskip2mm
\noindent
Let us now turn to the multidimensional variant of Theorems 1.1-1.3.
We will denote by $Z$ a standard normal random vector in $\R^d$, i.e., having
mean zero and an identity covariance matrix. Given i.i.d. random vectors 
$X,X_1,X_2,\dots$ in $\R^d$ with mean zero and identity covariance, consider 
the normalized sums
$$
Z_n = \frac{X_1 + \dots + X_n}{\sqrt{n}} \qquad (n = 1,2,\dots)
$$
We need to show that $D_\alpha(Z_n||Z)\to 0$ as $n\to\infty$, if and only if
$D_\alpha(Z_n||Z)$ is finite for some $n=n_0$, and
\be
\E\,e^{\left<X,t\right>} < e^{\beta |t|^2/2} \quad 
\text{for all} \ \ t \in \R^d, \ t\ne 0.
\en
Moreover, in this case $D_\alpha(Z_n||Z) = O(1/n)$, and $D_\alpha(Z_n||Z) = O(1/n^2)$
when the distribution of $X$ is symmetric about the origin. In fact,
a more precise Edgeworth-type expansion holds for $T_\alpha(Z_n||Z)$ 
in powers of $1/n$ similarly to (14.4)-(14.5), with the coefficients 
being polynomials of mixed cumulants of the components of $X$.

As for the proof of the theorems, much of the analysis developed before about 
the convergence in $T_\alpha$ (or $D_\alpha$), as well pointwise upper bounds on the
densities $p_n$ of $Z_n$, may easily be extended from dimension one to an arbitrary 
dimension $d$. Actually, the contractivity property of the functional $D_\alpha$ 
(Proposition 2.3) allows one to reduce the necessity part in Theorem 1.2 to 
the one dimensional case using a standard Wold type device. Indeed, consider 
the i.i.d. sequence $\left<X_i,\theta\right>$ with unit vectors $\theta$. 
Then, assuming that $D_\alpha(Z_n||Z)\to 0$ as $n \to \infty$, we get
$$
D_\alpha(\left<Z_n,\theta\right>||\left<Z,\theta\right>) \, \leq \, 
D_\alpha(Z_n||Z)\to 0.
$$
Since $\E \left<X_i,\theta\right> = 0$, $\E \left<X_i,\theta\right>^2 = 1$, and
$\left<Z,\theta\right> \sim N(0,1)$,
we are ready to apply the one dimensional variant of this theorem which gives
$$
\E\,e^{r\left<X,\theta\right>} < e^{\beta r^2/2} \quad \text{for all} \ \ r \ne 0.
$$
This is exactly the condition (16.1), thus proving the necessity part in Theorem 1.2.

Note that, as in dimension one (cf. Proposition 4.1), the finiteness of 
$D_\alpha(X||Z)$ guarantees that $\E\,e^{c|X|^2} < \infty$ for all $c < 1/(2\beta)$. 
In particular, the characteristic function $f(t) = \E\,e^{i\left<X,t\right>}$ 
now extends as an entire function to the $d$-dimensional complex space $\C^d$. 
Most important properties of the densities $p_n$ of $Z_n$ rely upon the function
$$
\psi(u) = f(iu) \,e^{-\beta |u|^2/2} = \E\,e^{-\left<X,u\right>}\,e^{-\beta |u|^2/2} 
\qquad (u \in \R^d).
$$

\vskip5mm
{\bf Lemma 16.1.} {\it If $T_\alpha = T_\alpha(X||Z) < \infty$, then $\psi(u)$ 
tends to zero as $|u| \rightarrow \infty$ and belongs to $L^k(\R^d)$ for any integer
$k \geq \alpha$. Moreover, up to some $(k,d)$-dependent constants $c_{k,d}$, we have
\be
\int_{\R^d} \psi(u)^k\,du \, \le \, c_{k,d}\, 
\big(1 + (\alpha - 1)T_\alpha\big)^{\frac{k}{\alpha}}.
\en
}

\vskip2mm
The first assertion is a multidimensional analog of Proposition 4.3; it can be
proved with very similar arguments as in dimension one. The second assertion 
generalizing Corollary 5.2 can be proved by using the contractivity properties 
of the $d$-dimensional Weierstrass transform
$$
W_t u(x) = \frac{1}{(2\pi t)^{d/2}} \int_{\R^d} e^{-\frac{|x-y|^2}{2t}}\,u(y)\,dy, 
\qquad x \in \R^d, \ t > 0.
$$
In particular, in $\R^d$ the inequality (5.1) takes the form
$
\E\,e^{\frac{1}{2\beta}\, |Z_k|^2} \leq c_{k,d}\, 
(1 + (\alpha - 1)T_\alpha)^{\frac{k}{\alpha}},
$
from which (16.2) easily follows. In case $\alpha = 2$, one may adapt Lemma 6.3 
as well to the multidimensional situation with its Parseval identity in $\R^d$. 
Furthermore, Proposition 6.2 is extended as
$$
\frac {1}{(2\pi)^{d/2}} \int_{\R^d} \psi(u)^2\,du \, \le \, 1+\chi^2(X,Z),
$$
thus refining (16.2) for $k=2$.

Repeating the arguments as in Section 13, one may also extend the corresponding
upper pointwise bounds on the densities.

\vskip5mm
{\bf Lemma 16.2.} {\it If $T_\alpha(X||Z) < \infty$, then
for all $x \in \R^d$ and $n \geq n_\beta = \max(\beta,2)$,
\be
p_n(x) \, \le \, A_{\alpha,d}\, n^{d/2}\, e^{-|x|^2/(2\beta)}\,
\psi\Big(-\frac{x}{\beta\sqrt n}\Big)^{n-n_\beta},
\en
where $A_{\alpha,d}$ depends on $(\alpha,d)$ only. In particular, there exist 
constants $x_0 > 0$ and $\delta \in (0,1)$ depending on the density $p$ such that
for all $n$ large enough
\be
p_n(x) \, \le \, 
\delta^n  e^{-|x|^2/(2\beta)}\,\psi\Big(-\frac{x}{\beta\sqrt n}\Big)^{n/2}, 
\quad 
whenever \ \ |x| \geq x_0 \sqrt{n}.
\en
}

\vskip2mm
{\bf Proof of Theorem 1.2} (Sufficiency part) {\bf and Theorem 1.3}. 
Assume that $n_0 = 1$. Hence $Z_n$ admits density $p_n$ for any $n \geq 1$. 
We need to derive the asymptotic behavior of
$$
(\alpha - 1)\, T_\alpha(Z_n||Z) = \int_{\R^d} w_n^\alpha(x)\, dx - 1, \qquad
w_n(x) = \frac{p_n(x)}{\varphi(x)^{1/\beta}},
$$
where $\varphi$ is the standard normal density on $\R^d$. To this aim, it is natural 
to split the integration into the four shell-type regions. The asymptotic behavior 
of the integrals
$$
I_0 = \int_{|x| < M_n} w_n^\alpha(x)\, dx, \qquad M_n = \sqrt{2(l-1)\log n},
$$
may be studied as in dimension one (cf. Lemma 11.2) by virtue of the Edgeworth 
expansion for $p_n(x)$ on the balls $|x| < M_n$ with a non-uniform error term.
To this aim, a multidimensional variant of Proposition 11.1 is used as stated 
in the monograph [BR-R], Theorem 19.2: Uniformly in $\R^d$
\be
p_n(x) = \varphi_s(x) + o\big(n^{-(s-2)/2}\big)\frac 1{1+|x|^s}, \qquad
\varphi_s(x) = \varphi(x) +  \varphi(x) \sum_{k=1}^{s-2}
\frac{q_k(x)}{n^{k/2}},
\en
where each $q_k$ represents a polynomial whose coefficients involve mixed cumulant
of the components of $X$ of order up to $k+2$. In particular, if the distribution 
of $X$ is symmetric about the origin, then $q_1(x) = 0$ and thus there is no 
$1/\sqrt{n}$ term in the sum (16.5).

In this way, we will arrive at the Edgeworth-type expansion for $I_0$ 
similarly to dimension one, which readily implies that $I_0 - 1 = O(1/n)$ in 
general, and $I_0 - 1 = O(1/n^2)$ when the distribution of $X$ is symmetric.
As a result, it remains to establish a polynomial smallness of the integrals
$$
I_1 = \int_{|x| > x_0 \sqrt{n}}\, w_n^\alpha(x)\, dx, \quad
I_2 = \int_{x_1 \sqrt{n} < |x| < x_0 \sqrt{n}}\, w_n^\alpha(x)\, dx, \quad
I_3 = \int_{M_n < |x| < x_1 \sqrt{n}}\, w_n^\alpha(x)\, dx
$$
with $x_1>0$ being any fixed small number, and $x_0>x_1$ depending on the density 
$p$. The bounds (16.2)-(16.4) allow us to properly estimate these integrals
as functions of $n$, by modifying the arguments from the previous section. 
Using (16.4) and (16.2) with $k = [\alpha] + 1$ and assuming that $\psi \leq 1$, 
we get for all $n$ large enough
$$
I_1 \, \leq \,
C_1 \int_{|x| > x_0 \sqrt{n}} 
\psi\Big(-\frac{x}{\beta\sqrt n}\Big)^{\alpha n/2}\, dx \, \leq \,
C_2\,\delta^{\alpha n}\, n^{d/2} \int_{\R^d} \psi(u)^k\, du \, \le \, C_3\,\delta_1^n
$$
with some constants $C_j,x_0>0$ and $0 < \delta < \delta_1 < 1$ which do not 
dependent on $n$.

For the region of $I_2$, thanks to condition (1.4), we have 
$\delta_2 = \max_{x_0 \leq |u| \leq x_1} \psi(u) < 1$. 
Hence, by (16.3), putting $n_1 = n - \max(\beta,2)$, we obtain that with 
some constants $C_j > 0$
\bee
I_2 
 & \le & 
C_1\, n^{d\alpha/2}
\int_{x_1 \sqrt{n} < |x| < x_0 \sqrt{n}} 
\psi\Big(-\frac{x}{\beta\sqrt n}\Big)^{n_1}\,dx \\
 & = &
C_2\, n^{d(\alpha + 1)/2} 
\int_{x_1 \sqrt{n} < |x| < x_0 \sqrt{n}} \psi(u)^{n_1}\,du
 \ \le \ 
C_3\,n^{d(\alpha + 2)/2}\,x_0^d\,\delta_2^{n_1}
\ene
which is decaying exponentially fast like $I_1$. 

Finally, using the analyticity of $f$, we have
$\psi(u) \leq e^{-(\beta - 1) |u|^2/4}$ in a sufficiently small ball $|u|<r$, 
so that
$$
\psi\Big(-\frac{x}{\beta\sqrt n}\Big)^{n_1} \leq 
\psi\Big(-\frac{x}{\beta\sqrt n}\Big)^{n/2} \leq 
e^{-x^2/(8 \alpha \beta)}, \qquad |x| < \beta r\sqrt n,
$$
for all $n \geq 2\max(\beta,2)$. Therefore, by (16.2), in this ball
$
w_n(x) \le A_{\alpha,d}^\alpha\, n^{d\alpha/2}\, e^{-|x|^2/(8 \alpha \beta)},
$
which gives with $x_1 = \beta r$
\bee
I_3 
 & \le &
C_1\, n^{d\alpha/2} \int_{M_n < |x| < x_1\sqrt n}
e^{-|x|^2/(8 \alpha \beta)}\,dx \\
 & < &
C_2\, n^{d\alpha/2}\,\P\{|Z|^2 > M_n^2/(8d\alpha \beta)\} \, \leq \,
C_3\, n^{d\alpha/2}\, e^{-M_n^2/(8d \alpha \beta)}
\, = \, 
C_3\, n^{-(\frac{l - 1}{4d\alpha \beta} - \frac{\alpha}{2})}.
\ene
Collecting these bounds, we get that $I_1 + I_2 + I_3 = o(n^{-l/8 d \alpha \beta})$
for all sufficiently large $l$, thus proving Theorem 1.2. 

For the proof of Theorem 1.3 in $\R^d$, we need to investigate the suprema
$$
J_0 = \sup_{|x| \le M_n} \, \frac{|p_n(x) - \varphi_s(x)|}{\varphi(x)^{1/\beta}}, 
\qquad
J_1  = \sup_{|x| \geq x_1\sqrt n} \, \frac{p_n(x)}{\varphi(x)^{1/\beta}}, \qquad
J_2 = \sup_{M_n \leq |x| \leq x_1\sqrt n}\,  \frac{p_n(x)}{\varphi(x)^{1/\beta}}
$$ 
with some $x_1 > 0$ and assuming that $M_n < x_1 \sqrt n$. An application of the
expansion (16.5) implies that $J_0 = O(1/\sqrt{n})$ in general and
$I_0 = O(1/n)$ when the distribution of $X$ is symmetric.
The polynomial smallness of $J_1$ and $J_2$ (for sufficiently large values
of $l$ in the definition of $M_n$) follows from Lemma 16.2, by repeating 
the arguments of the proof of Theorem 15.1.
\qed

\vskip10mm
\section{{\bf Some Examples and Counter-Examples}}
\setcounter{equation}{0}

\vskip2mm
\noindent
Given a random variable $X$ such that $\E X=0$, $\E X^2=1$, consider the function 
$\psi(t)=e^{-t^2}\, \E\, e^{tX}$ ($t \in \R$). As before, put
$$
Z_n = \frac{X_1 + \dots + X_n}{\sqrt{n}},
$$
where $X_j$'s are independent copies of $X$. One immediate consequence of 
Theorem 1.1 (with $n_0 = 1$) is the following characterization.

\vskip5mm
{\bf Theorem 17.1.} {\it Assume that the random variable $X$ has a density $p$
such that
\be
\int_{-\infty}^\infty p(x)^2\,e^{x^2/2}\,dx < \infty.
\en
Then $\chi^2(Z_n,Z) \rightarrow 0$ as $n \rightarrow \infty$ for $Z \sim N(0,1)$, 
if and only if
\be
\psi(t) < 1 \quad for \ all \ \ t \neq 0.
\en
}

The assumption (17.1) is fulfilled, for example, when $X$ is bounded and has a square 
integ\-rab\-le density. We now illustrate Theorem 17.1 and the more general Theorem 1.2 
with a few examples (mostly in dimension one).

\vskip2mm
{\bf Uniform distribution.} Let $X$ be uniformly distributed on the segment 
$[-\sqrt 3,\sqrt 3]$. The characteristic function of $X$ is given by 
$f(t)=\sin(t\sqrt 3)/(t\sqrt 3)$, and for imaginary
values $t=iy$, we have the simple estimate
\be
f(iy)=\frac{\sinh(y\sqrt 3)}{y\sqrt 3} < e^{y^2/2},\qquad y\in\mathbb R \ (y \neq 0),
\en
so that (17.2) does hold. In this case the first moments are given by 
$\alpha_2 = 1$, $\alpha_3 = 0$, $\alpha_4 = \frac{9}{5}$.
Therefore, by Theorem 17.1, $\chi^2(Z_n,Z) \rightarrow 0$ as $n \rightarrow \infty$.
Moreover, Theorem 1.1 provides an asymptotic expansion (1.3) which becomes
$$
\chi^2(Z_n,Z) \, = \, \frac{3}{50\,n^2} + O\Big(\frac{1}{n^3}\Big).
$$

In fact, the property (17.3) means that the condition (1.4) of a more general Theorem 1.2 
is fulfilled in the whole range of indexes $\alpha > 1$. Using the formula (14.5), we 
therefore obtain a stronger assertion 
$T_\alpha(Z_n||Z) = \frac{\alpha}{2}\,\chi^2(Z_n,Z) + O(\frac{1}{n^3})$,
and a similar one for $D_\alpha$.

\vskip2mm
{\bf Convex mixtures of centered Gaussian measures.} 
Consider the densities of the form
$$
p(x) = \int_0^\infty \frac{1}{\sigma \sqrt{2\pi}}\,e^{-x^2/2\sigma^2}
\,d\pi(\sigma^2),  \qquad x \in \R,                                     
$$
where $\pi$ is a (mixing) probability measure on the positive half-axis with
$
\int_0^\infty \sigma^2 d\pi(\sigma^2) = 1.                             
$
The random variable $X$ with this density has mean zero and variance one, and
its distribution is equal to that of $\sqrt{\xi}\,Z$, where $\xi$ is independent of
$Z \sim N(0,1)$ and is distributed according to $\pi$. As in Example 9.3, 
$\chi^2(Z_n,Z) < \infty$ for some $n = n_0$, if and only if 
$\pi$ is supported on the interval $(0,2)$, and its distribution function 
$F(\ep) = \pi((0,\ep])$ satisfies
\be
\inf_n \int_0^1 \frac{F(\ep)^{2n}}{\ep^{3/2}}\,d\ep < \infty, \qquad
\inf_n \int_1^2 \frac{(1 - F(\ep))^{2n}}{(2 - \ep)^{3/2}}\,d\ep < \infty.
\en

On the other hand, the distribution of $X$ has the Laplace transform
$$
\E\,e^{tX} = \int_0^{\infty} e^{\sigma^2 t^2/2} \,d\pi(\sigma^2) = 
\E\,e^{\xi t^2/2},  \qquad t \in \R.
$$
Hence, the condition $\chi^2(Z_n,Z) < \infty$ guarantees that (17.2) is fulfilled.
Without that condition,
$\E\,e^{tX} < e^{t^2}$ for all $t \neq 0$, if and only if $\P\{\xi \leq 2\} = 1$
and $\P\{\xi = 2\} < 1$. Here, $\P\{\xi = 2\} = 1$ is not possible in view of the
second moment assumption $\E X^2 = \E\,\xi = 1$. 

Hence, one concludes that
$\chi^2(Z_n,Z) \rightarrow 0$ as $n \rightarrow \infty$, if and only if the measure
$\pi$ is supported on the interval $(0,2)$ and satisfies the condition (17.4). 
In this case, we obtain the expansion (1.3) which reads
$$
\chi^2(Z_n,Z) \, = \, \frac{3\,(m-1)^2}{8\,n^2} + O\Big(\frac{1}{n^3}\Big), \qquad
m = \int_0^\infty \sigma^4 \,d\pi(\sigma^2).
$$

\vskip2mm
{\bf Distributions with Gaussian component.}
Consider random variables of the form
$$
X = a\xi + bZ \qquad (a^2 + b^2 = 1, \ a,b>0)
$$
assuming that $\E \xi = 0$, $\E\xi^2 = 1$, and where 
$Z \sim N(0,1)$ is independent of $\xi$. The distribution of $X$ is a convex 
mixture of shifted Gaussian measures on the line with variance $b^2$. It admits a density
$$
p(x) = \frac{1}{b}\,\E\,\varphi\Big(\frac{x - a\xi}{b}\Big), \qquad x \in \R.
$$
To ensure finiteness of $\chi^2(X,Z)$ (and even finiteness of $\chi^2(Z_n,Z)$ with 
some $n$), the random variable $\xi$ should have a finite Gaussian moment, or equivalently,
the Laplace transform of the distribution of $\xi$ should admit a subgaussian bound
\be
\E\,e^{t \xi} \leq e^{\sigma^2 t^2/2}, \qquad t \in \R,
\en
with some finite $\sigma>0$. Let $\sigma$ be an optimal value in this inequality 
(necessarily $\sigma \geq 1$). It then follows that
$\E\,e^{c \xi^2} < \infty$ whenever $c < 1/(2\sigma^2)$.

Squaring the formula for $p(x)$, we easily find an expression for the $\chi^2$-distance,
namely,
$$
1 + \chi^2(X,Z) \, = \, \frac{1}{\sqrt{1 - a^4}}\,\E\,\exp\bigg\{\frac{a^2}{2(1 - a^2)}\,
\Big(\frac{2}{1 + a^2}\,(\xi + \eta)^2 - (\xi^2 + \eta^2)\Big)\bigg\},
$$
where $\eta$ is an independent copy of $\xi$. Using $(\xi + \eta)^2 \leq 2\xi^2 + 2\eta^2$,
we are lead to a simple upper bound
$$
1 + \chi^2(X,Z) \leq 
\frac{1}{\sqrt{1 - a^4}}\,\Big(\E\,e^{\frac{a^2}{2(1 + a^2)}\,\xi^2}\Big)^2.
$$
Hence, $\chi^2(X,Z) < \infty$ whenever $a < a_\sigma = \frac{1}{\sqrt{\sigma^2 - 1}}$,
which is automatically fulfilled in case $\sigma^2 \leq 2$. Moreover, for all $t \neq 0$,
$$
\E\,e^{tX} = \E\,e^{at\xi}\,e^{b^2t^2/2} \leq e^{(\sigma^2 a^2 + b^2)\, t^2/2} 
= e^{((\sigma^2 - 1) a^2 - 1)\, t^2/2} < e^{t^2}
$$
under the same constraint $a < a_\sigma$. Thus we conclude, by applying Theorem 17.1, that
$\chi^2(Z_n,Z) \rightarrow 0$ as $n \rightarrow \infty$, if $a < \frac{1}{\sqrt{\sigma^2 - 1}}$.
In case $\sigma^2 \leq 2$, this convergence holds for all admissible parameters $(a,b)$.

\vskip2mm
{\bf Distributions with finite Gaussian moment.} Suppose that a random variable $X$
with mean zero and variance one has finite Gaussian moment $M = \E\,e^{cX^2}$ ($c>0$).
It is well-known that the property (17.5) is fulfilled for some $\sigma \geq 1$; moreover,
one can show that an optimal value satisfies $\sigma^2 \leq \frac{4\log M}{c\log 2}$.
This means that condition (1.4) is fulfilled for any $\alpha>1$ such that
$\beta < \sigma^2$. We conclude that, if $D_\alpha(X||Z) < \infty$, then
$D_\alpha(Z_n||Z) \rightarrow 0$ with any $\alpha < \frac{\sigma^2}{\sigma^2 - 1}$.

\vskip2mm
{\bf Conditions in terms of exponential series.}
Consider a symmetric density of the form
$$
p(x) = \varphi(x) \sum_{k=0}^\infty \frac{\sigma_k}{2^k k!} H_{2k}(x),
\qquad x\in\R,
$$
with $\sigma_0 = 1$ and $\sigma_1=0$ (which means that $\E X^2 = 1$ for
the random variable with density $p$).
In view of Section 6, condition (17.1) is fulfilled, if and only if the
series
$$
\chi^2(X,Z) \, = \, \sum_{k=2}^\infty \frac{(2k)!}{4^k\,k!^2}\,\sigma_k^2 
\, \sim \, \sum_{k=2}^\infty \frac{1}{\sqrt{k}}\,\sigma_k^2
$$
is convergent (which is fulfilled automatically, when $p$ is compactly supported
and bounded). Assuming additionally that $\sup_{k\ge 2} \sigma_k \le 1$, we also have
$$
\E\,e^{tX} \, = \, e^{t^2/2}\,
\bigg[1+\sum_{k=2}^{\infty}\frac{\sigma_k}{k!}\Big(\frac{t^2}2\Big)^k\bigg]
 \, \le \, e^{t^2/2}\Big(e^{t^2/2}-\frac{t^2}{2}\Big) \, < \, e^{t^2}, \qquad t\ne 0.
$$
Hence, in this case, by Theorem~17.1, $\chi^2(Z_n,Z) \rightarrow 0$ as $n \rightarrow \infty$.
Moreover,  according to the expansion (1.3), we have $\chi^2(Z_n,Z) = O(1/n^2)$.
This assertion strengthens the result of [F] (under weaker assumptions).

\vskip2mm
{\bf Log-concave probability distributions.} More examples including those in higher
dimensions illustrate the multidimensional Theorem 1.2
within the class of densities $p(x) = e^{-V(x)}$ supported on 
some open convex region $\Omega \subset \R^d$, where $V$ is a $C^2$-convex 
function with Hessian satisfying $V''(x) \geq c\, {\rm I}_d$ in the sense of positive
definite matrices ($c > 0$). The probability measures with such
densities are known to admit logarithmic Sobolev inequalities (via the Bakry-Emery criterion).
In particular, they satisfy transport-entropy inequalities which in turn
can be used to get a subgaussian bound on the Laplace transform such as
$$
\E\,e^{r u(X)} \leq e^{r^2/(2c)}, \qquad r \in \R.
$$
Here, $u$ may be an arbitrary function on $\R^d$ with Lipschitz semi-norm
$\|u\|_{\rm Lip} \leq 1$, such that $\E\,u(X) = 0$ (cf. [B-G], [O-V]).
In particular, if $\E X = 0$, one may choose an arbitrary linear function 
$u(x) = \left<x,\theta\right>$ with $|\theta| = 1$. Hence, the condition (1.4)  will be 
fulfilled, as long as $c > \frac{1}{\beta}$. Moreover, the property $D_\alpha(X||Z) < \infty$
will also hold in this case, since necessarily
$$
V(x) \geq V(x_0) + \left<V'(x_0), x-x_0\right> + \frac{c}{2}\,|x-x_0|^2
$$
for all $x,x_0 \in \Omega$. Applying Theorem 1.2, we get:

\vskip5mm
{\bf Corollary 17.2.} {\it If a random vector $X$ in $\R^d$ with mean zero and identity
covariance matrix has density $p = e^{-V}$ such that $V'' \geq c\, {\rm I}_d$ 
$(0 < c \leq 1)$ on the supporting open convex region,
then $D_\alpha(Z_n||Z) \rightarrow 0$ as $n \rightarrow \infty$, whenever
$\alpha < \frac{1}{1-c}$.
}

\vskip10mm
\section{{\bf Convolution of Bernoulli with Gaussian}}
\setcounter{equation}{0}

\vskip2mm
\noindent
One might wonder whether or not it is possible to replace the condition (1.1)
in Theorem 1.1 with a slightly weaker requirement like $\E\,e^{tX} \leq e^{t^2}$
(hoping e.g. that the strict inequality would automatically hold, in view of the
assumption $\E X^2 = 1$). The answer is negative, as the following statement
shows:

\vskip5mm
{\bf Proposition 18.1.} {\it There exists a random variable $X$ with $\E X = 0$, $\E X^2 = 1$,
$\chi^2(X,Z) < \infty$ for $Z \sim N(0,1)$, and such that the inequality
\be
\E\,e^{tX} < e^{t^2}
\en
is fulfilled for all $t \neq 0$ except for exactly one point $t_0 \neq 0$.
}

\vskip5mm
Since (18.1) is violated (although at one point only), Theorem 1.1 implies that
convergence $\chi^2(Z_n,Z) \rightarrow 0$ does not hold any more.

Let us describe explicitly one family of distributions satisfying the assertion of this proposition.
Returning to one of the previous examples, consider random variables of the form
$$
X_p = a\xi + bZ \qquad (a,b>0),
$$
assuming that $\xi$ takes two values $q$ and $-p$ with probabilities 
$p$ and $q$, respectively ($p,q > 0$, $p + q = 1$), and where $Z \sim N(0,1)$ is independent 
of $\xi$. Clearly, $\E X_p = 0$, and we have the constraint
\be
\E X_p^2 = pq\, a^2 + b^2 = 1.
\en
The density $w$ of $X_p$ represents a convex mixture of two shifted Gaussian densities,
$$
w(x) = \frac{p}{b}\,\varphi\Big(\frac{x - aq}{b}\Big) + 
\frac{q}{b}\,\varphi\Big(\frac{x + ap}{b}\Big),
$$
and the condition $\chi^2(X,Z) < \infty$ obviously holds (since necessarily $b<1$). 

Now, let $\sigma^2 = \sigma^2(p,q)$ denote the smallest positive
constant such that the following inequality holds
\be
\E\,e^{t\xi} = p e^{qt} + q e^{-pt} \leq e^{\sigma^2 t^2/2}, \qquad t \in \R.
\en
This is the so-called subgaussian constant for the Bernoulli distribution.
Since $\E\,e^{tX_p} = \E\,e^{at\xi}\, e^{b^2 t^2/2}$, (18.3) yields
$$
\E\,e^{tX_p} \leq e^{(\sigma^2 a^2 + b^2)\, t^2/2}, \qquad t \in \R,
$$
with an optimal constant $\sigma^2 a^2 + b^2$ in the exponent on the right-hand side.
Thus, according to the requirement (18.1), we get another constraint 
$\sigma^2 a^2 + b^2 = 2$. Combining it with (18.2), we find that necessarily
$$
a^2 = \frac{1}{\sigma^2 - pq}, \qquad
b^2 = \frac{\sigma^2 - 2pq}{\sigma^2 - pq},
$$
which makes sense provided that $\sigma^2 > 2pq$. The subgaussian constant 
for the Bernoulli distribution is known to be (cf. [B-H-T], Proposition 2.3)
$$
\sigma^2 = \frac{p-q}{2\,(\log p - \log q)}.
$$
It is easy to see that (18.3) becomes equality for 
$t_0 = -2\,(\log p - \log q)$, which is a unique non-zero point with such property, 
as long as $p \neq q$.

Hence we conclude that the random variable $X = X_p$ satisfies the assertion of 
Proposition 18.1, if and only if 
\be
\frac{p-q}{2\,(\log p - \log q)} > 2pq.
\en
This inequality holds, provided that $p$ is sufficiently close to 0 or 1
(although it is not true for a neighborhood of $1/2$). For example, one may choose
$p = 1/6$. More precisely, for some constant $p_0 \in (0,\frac{1}{2})$,
(18.4) holds for $p$ from the set $(0,p_0) \cup (1-p_0,1)$, while for $p$ from 
$(p_0,1-p_0)$ it holds with an opposite inequality sign.

\end{document}